\documentclass{article}
\usepackage{amsmath,amssymb,latexsym,a4}
\usepackage[title]{appendix}

\newtheorem{theorem}{Theorem}
\newtheorem{lemma}{Lemma}[section]

\newtheorem{corollary}[lemma]{Corollary}
\newtheorem{proposition}[lemma]{Proposition}

\newtheorem{definition}{Definition}
\newtheorem{example}[lemma]{Example}
\newtheorem{remark}[lemma]{Remark}
\newtheorem{conjecture}{Conjecture}

\newcommand{\bl}{\begin{lemma}}
\newcommand{\el}{\end{lemma}}
\newcommand{\bt}{\begin{theorem}}
\newcommand{\et}{\end{theorem}}
\newcommand{\bcor}{\begin{corollary}}
\newcommand{\ecor}{\end{corollary}}
\newcommand{\bp}{\proof{.}}
\newcommand{\ep}{\eop}
\newcommand{\bpr}{\begin{proposition}}
\newcommand{\epr}{\end{proposition}}
\newcommand{\brem}{\begin{remark} \em}
\newcommand{\erem}{\end{remark}}
\newcommand{\bd}{\begin{definition} \em}
\newcommand{\ed}{\end{definition}}
\newcommand{\bex}{\begin{example} \em
}
\newcommand{\eex}{\end{example}}
\newcommand{\beq}{\begin{equation} }
\newcommand{\eeq}{\end{equation}}

\newcommand{\bi}{\begin{itemize}
  }
\newcommand{\ei}{\end{itemize}}
\newcommand{\ben}{\begin{enumerate} }
\newcommand{\een}{\end{enumerate} }

\newcommand{\refeq}[1]{(\ref{#1})}
\newenvironment{enumr}{

\begin{enumerate}     }{\end{enumerate}

}
\newenvironment{enuma}{

\begin{enumerate}     }{\end{enumerate}

}

\newcommand{\benr}{\begin{enumr}
  }
\newcommand{\eenr}{
\end{enumr}}

\newcommand{\ignore}[1]{}

\newcommand{\al}[1]{\forall #1\:}
\newcommand{\ex}[1]{\exists #1\:}

\newlength{\hilflh}

\newcommand{\naturals}{\mathbb{N}}

\newcommand{\bb}{\mbox{``}}
\newcommand{\qq}{\mbox{''}}
\newcommand{\cL}{{\mathcal L}}

\newcommand{\cE}{{\mathcal E}}

\newcommand{\cP}{{\mathcal P}}

\newcommand{\cW}{{\mathcal W}}

\newcommand{\cI}{\mathcal{I}}

\newcommand{\fB}{\mathfrak{B}}
\newcommand{\fL}{\mathfrak{L}}
\newcommand{\fG}{\mathfrak{G}}
\newcommand{\fI}{\mathfrak{I}}
\newcommand{\fT}{\mathfrak{T}}
\newcommand{\fC}{\mathfrak{C}}

\newcommand{\fO}{\mathfrak{O}}
\newcommand{\fM}{\mathfrak{M}}

\newcommand{\oI}{\bar{\mI}}

\newcommand{\ga}{\alpha}
\newcommand{\gb}{\beta}
\newcommand{\gd}{\delta}
\renewcommand{\ge}{\varepsilon}
\newcommand{\gl}{\lambda}

\newcommand{\gs}{\sigma}
\newcommand{\gy}{\gamma}
\newcommand{\gw}{\omega}
\newcommand{\gS}{\Sigma}

\renewcommand{\phi}{\varphi}

\newcommand{\eqv}{\leftrightarrow}

\newcommand{\thr}[1]{[#1]}

\newcommand{\ol}{\overline}

\newcommand{\Glp}{\mathrm{GLP}}

\newcommand{\Con}{\mathrm{Con}}

\newcommand{\Prf}{\mathrm{Prf}}
\newcommand{\Ax}{\mathrm{Ax}}

\newcommand{\PRA}{\mathsf{PRA}}

\newcommand{\PA}{\mathsf{PA}}
\newcommand{\EA}{\mathrm{EA}}

\newcommand{\gn}[1]{\ulcorner #1 \urcorner}

\newcommand{\ord}[2]{|#2|_{\Pi_1^0}}

\newcommand{\fc}{\Vdash}      
\renewcommand{\models}{\vDash}      
\newcommand{\nfc}{\nVdash}

\newcommand{\Imp}{\Rightarrow}

\newcommand{\Var}{\mathrm{Var}}

\newcommand{\Lim}{\mathrm{Lim}}
\newcommand{\Suc}{\mathrm{Suc}}

\newcommand{\iffdef}{\stackrel{\text{def}}{\iff}}

\newcommand{\nat}{\naturals}

\renewcommand{\leq}{\leqslant}
\renewcommand{\geq}{\geqslant}

\newcommand{\Rc}{\mathrm{RC}}

\newcommand{\Ig}{\cI}
\newcommand{\oIg}{\overline{\Ig}}

\newcommand{\TRC}{\fT^0_{\Rc}}
\newcommand{\Sat}{\mathrm{Sat}}
\newcommand{\mO}{\mathrm{O}}
\newcommand{\mI}{\mathrm{I}}
\newcommand{\spec}{\mathrm{sp}}
\newcommand{\thh}{\mathrm{th}}

\newcommand{\eop}{$\Box$ \protect\par \addvspace{\topsep}}
\newcommand{\proof}[1]{\protect\par\addvspace{\topsep}\noindent {\bf Proof#1}}

\newcommand{\Rcn}{{\mathrm{RC}^\nab}}
\newcommand{\ofG}{\ol{\fG}}
\newcommand{\iPi}{\mathrm{\Pi}}
\newcommand{\Kp}{\mathrm{K}^+}
\newcommand{\Kxp}{\mathrm{K4}^+}
\newcommand{\nab}{\nabla\hspace{-0.7pt}}
\newcommand{\Wo}{\mathbb{W}}
\newcommand{\Fo}{\mathbb{F}}
\newcommand{\Fon}{\mathbb{F}^\nab}
\newcommand{\mR}{\mathrm{R}}
\newcommand{\mQ}{\mathrm{Q}}
\newcommand{\mA}{\mathrm{A}}
\newcommand{\mC}{\mathrm{C}}
\newcommand{\BS}{\mathrm{B\Sigma}}

\begin{document}

\title{Reflection calculus and conservativity spectra}

\author{Lev D. Beklemishev\thanks{Research financed by a grant of the Russian Science Foundation (project No. 16-11-10252).} \\
Steklov Mathematical Institute of
Russian Academy of Sciences \\ Gubkina str. 8, Moscow,\ \texttt{bekl@mi.ras.ru}}

\maketitle

\begin{abstract}
Strictly positive logics recently attracted attention both in the description logic and in the provability logic communities for their combination of efficiency and sufficient expressivity. The language of Reflection Calculus RC consists of implications between formulas built up from propositional variables and constant `true' using only conjunction and diamond modalities which are interpreted in Peano arithmetic as restricted uniform reflection principles.

We extend the language of $\Rc$ by another series of modalities representing the operators associating with a given arithmetical theory $T$ its fragment axiomatized by all theorems of $T$ of arithmetical complexity $\Pi^0_n$, for all $n>0$. We note that such operators, in a strong sense, cannot be represented in the full language of modal logic.

We formulate a formal system $\Rcn$ extending $\Rc$ that is sound and, as we conjecture, complete under this interpretation. We show that in this system one is able to express iterations of reflection principles up to any ordinal $<\ge_0$. Secondly, we provide normal forms for its variable-free fragment. Thereby, this fragment is shown to be algorithmically decidable and complete w.r.t.\ its natural arithmetical semantics.

In the last part of the paper we characterize in several natural ways the Lindenbaum--Tarski algebra of the variable-free fragment of $\Rcn$ and its dual Kripke structure. Most importantly, the elements of this algebra correspond to the sequences of proof-theoretic $\Pi^0_{n+1}$-ordinals of bounded fragments of Peano arithmetic called \emph{conservativity spectra}, as well as to the points of the well-known Ignatiev Kripke model.
\end{abstract}

\section{Introduction}

A system, called \emph{Reflection Calculus} and denoted
$\Rc$, was introduced in \cite{Bek12a} and, in a slightly different format, in \cite{Das12}.
From the point of view of modal logic, $\Rc$ can be seen as a
fragment of Japaridze's polymodal provability logic $\Glp$ \cite{Dzh86,Boo93} consisting of the implications of the form $A\to B$, where $A$ and $B$ are formulas built-up from $\top$ and propositional variables using just $\land$ and the diamond modalities. We call such formulas $A$ and $B$ \emph{strictly positive}.

Strictly positive modal logics, earlier and in a different guise, appeared in the work on description logic. They serve as a good compromise between the concerns of efficiency and sufficient expressivity in the knowledge base query answering. In particular, the strictly positive language corresponds to the OWL2EL profile of the OWL web ontology language, and is used in large ontology bases such as SNOMED CT. The papers \cite{KKTWZ,KKTWZ-arx} undertake a general study of strictly positive logics and provide more references, especially in the description logic and in the universal algebraic traditions.

Our concerns in the development of strictly positive provability logic are, in a sense, similar. Reflection calculus $\Rc$ is much simpler than its modal companion $\Glp$ yet expressive enough for its main proof-theoretic applications. It has been outlined in \cite{Bek12a} that $\Rc$
allows to define a natural system of ordinal notations up to
$\ge_0$ and serves as a convenient basis for a proof-theoretic
analysis of Peano Arithmetic in the style of \cite{Bek04,Bek05en}.
This includes a consistency proof for Peano arithmetic based on transfinite induction up to $\ge_0$, a characterization of its
$\Pi_n^0$-consequences in terms of iterated reflection principles, a slowly terminating term rewriting system~\cite{BekOn15} and a combinatorial independence result~\cite{Bek06}.

An axiomatization of $\Rc$ (as an equational calculus) has been found by Evgeny Dashkov in his paper~\cite{Das12} which initiated the study of strictly positive provability logics. Dashkov proved two important further facts about $\Rc$ which sharply contrast with the
corresponding properties of $\Glp$. Firstly, $\Rc$ is complete with respect to a natural class of finite Kripke frames. Secondly, $\Rc$ is decidable in polynomial time, whereas most of the standard modal logics are \textsc{PSpace}-complete and the same holds for the variable-free fragment of GLP \cite{Pakh14}.

Another advantage of going to a strictly positive language is
exploited in the present paper. Strictly positive modal formulas
allow for more general arithmetical interpretations than those of
the standard modal logic language. In particular, propositional
formulas can now be interpreted as arithmetical \emph{theories}
rather than individual \emph{sentences}. (Notice that the `negation' of a theory would not be well-defined.) As the first meaningful example for this framework we analysed an extension of $\Rc$ by a modality representing the full arithmetical uniform reflection principle \cite{Bek14}. The corresponding strictly positive logic, though arithmetically complete, complete w.r.t.~a nice class of finite Kripke models and polytime decidable, turned out not to be equivalent to the fragment of any standard normal modal logic.\footnote{This has not been noted in~\cite{Bek14}, however it follows from Theorem 3 of~\cite{Bek18b} saying that a s.p.\ logic is a fragment of a normal modal logic iff it is Kripke frame complete. Modulo some reformulations this result is, in fact, equivalent to Theorem 1 of~\cite{KKTWZ}.}

More generally, any monotone operator acting on the semilattice of arithmetical theories can be considered as a modality in strictly positive logic. One such operation is particularly attractive from the point of view of proof-theoretic applications, namely the map associating with a theory $T$ its fragment $\iPi_{n+1}(T)$ axiomatized by all theorems of $T$ of arithmetical complexity $\Pi^0_{n+1}$. Since the $\Pi^0_{n+1}$-conservativity relation of $T$ over $S$ can be expressed by $S\vdash \iPi_{n+1}(T)$, we call such operators \emph{$\Pi^0_{n+1}$-conservativity operators}.

This relates our study to the fruitful tradition of research on conservativity and interpretability logics, see e.g.\ \cite{Vis97,Vis08b,DJ,JooG11,HMa,HMb,Ign91}. Our framework happens to be both weaker and stronger than the traditional one: in our system we are able to express the conservativity relations for each class $\Pi^0_{n+1}$ and are able to relate not only sentences but theories. However, in this framework the negation is lacking and the conservativity is not a binary modality and cannot be iterated. Yet, we believe that the strictly positive language is both simpler and better tuned to the needs of proof-theoretic analysis of formal systems of arithmetic.

We introduce the system $\Rcn$ with modalities $\Diamond_n$ representing uniform reflection principles of arithmetical complexity $\Sigma_n$, and $\nab_n$ representing $\Pi_{n+1}$-conservati\-vi\-ty operators. We provide an adequate semantics of $\Rcn$ in terms of the semilattice $\fG_\EA$ of (numerated) arithmetical r.e.\ theories extending elementary arithmetic $\EA$. Further, we introduce transfinite iterations of monotone semi-idempotent operators along elementary well-orderings, somewhat generalizing the notion of a Turing--Feferman recursive progression of axiomatic systems but mainly following the same development as in \cite{Bek06}. Our first result shows that  $\ga$-iterations of modalities $\Diamond_n$, for each $n<\gw$ and ordinals $\ga<\ge_0$, are expressible in the algebra $\fG_\EA$. A variable-free strictly positive logic where such iterations are explicitly present in the language has been introduced by Hermo Reyes and Joosten \cite{JooRey16} which is, thereby, contained in $\Rcn$. However, possible generalisations of their system to larger ordinal notation systems would be out of scope of $\Rcn$.\footnote{In the latest version of their paper Hermo Reyes and Joosten did, in fact, exted their setup to arbitrary ordinal notation systems.}

Then we turn to a purely syntactic study of the variable-free fragment of the system $\Rcn$ and provide unique normal forms for its formulas. A corollary is that the variable-free fragment of $\Rcn$ is decidable and arithmetically complete.

Whereas the normal forms for the variable-free formulas of $\Rc$ correspond in a unique way to ordinals below $\ge_0$, the normal forms of $\Rcn$ are more general. It turns out that they are related in a canonical way to the collections of proof-theoretic ordinals of (bounded) arithmetical theories for each complexity level $\Pi_{n+1}$, as defined in \cite{Bek06}.

Studying the collections of proof-theoretic ordinals corresponding to several levels of logical complexity as single objects seems to be a rather recent and interesting development. Such collections appeared for the first time in the work of  Joost Joosten \cite{Joo15a} under the name \emph{Turing--Taylor expansions}. He established a one-to-one correspondence between such collections (for a certain class of theories) and the points of the Ignatiev universal model for the variable-free fragment of $\Glp$. We call such collections \emph{conservativity spectra} of arithmetical theories. Our results show that $\Rcn$ provides a way to syntactically represent and conveniently handle such conservativity spectra.

The third part of our paper provides an algebraic model $\fI$ for the variable-free fragment of $\Rcn$. This model is obtained in a canonical way on the basis of the Ignatiev model. Our main theorem states the isomorphism of several representations of $\fI$: the Lindenbaum--Tarski algebra of the variable-free fragment of $\Rcn$; a constructive representation in terms of sequences of ordinals below $\ge_0$; a representation in terms of the semilattice of bounded $\Rc$-theories and as the algebra of cones of the Ignatiev model. 
In Section \ref{ukf} we consider its dual relational structure $\fI^*$, which is universal for the variable-free fragment of $\Rcn$. We give a constructive characterization of this large Kripke frame in terms of sequences of ordinals.

Parts of this paper previously appeared in conference proceedings~\cite{Bek17c,Bek18a} though underwent a thorough revision here.
Thanks are due to Albert Visser for suggesting many improvements including Lemma \ref{vis}, as well as to Ilya Shapirovsky, Joost Joosten, and Evgeny Kolmakov for comments and corrections.

\section{The lattice of arithmetical theories}

We define the intended arithmetical interpretation of the strictly positive
modal language. Propositional variables (and strictly
positive formulas) will now denote possibly infinite theories rather than
individual sentences. We deal with r.e.\ theories formulated in the language of \emph{elementary arithmetic $\EA$} and containing the axioms of $\EA$. The theory $\EA$, aka $I\Delta_0(\exp)$ or $\text{EFA}$, is formulated in the language of Peano arithemtic enriched by a symbol for exponentiation ($2^x$). In addition to the standard quantifier-free defining axioms for all the symbols of the language, it contains the induction schema for bounded formulas (cf \cite{HP,Bek05en}). Bounded formulas in the language of $\EA$ are called \emph{elementary formulas}, the class of all such formulas is usually denoted $\Delta_0(\exp)$.

To avoid well-known problems with the representation of theories in arithmetic, we assume that each theory $S$ comes equipped with an
\emph{elementary numeration}, that is, a bounded formula $\gs(x)$ in the language of $\EA$ defining the set of axioms of $S$ in the standard model of arithmetic $\nat$.

Given such a $\gs$, we have a standard arithmetical $\Sigma^0_1$-formula $\Box_\gs(x)$ expressing the provability of $x$ in $S$ (see \cite{Fef60}). We often write
$\Box_\gs\phi$ for $\Box_\gs(\gn{\phi})$. The expression $\bar n$
denotes the numeral ${0'}^{\ldots}{'}$ ($n$ times). If $\phi(v)$ contains a
parameter $v$, then $\Box_\gs\phi(\bar x)$ denotes a formula (with a
parameter $x$) expressing the provability of the sentence $\phi(\bar x/v)$ in $S$.

Given two numerations $\gs$ and $\tau$, we write $\gs\leq_{\EA}
\tau$ if $$\EA\vdash \al{x}(\Box_\tau(x)\to \Box_\gs(x)).$$
We will only consider the numerations
$\gs$ such that $\gs\leq_\EA \gs_\EA,$ where $\gs_\EA$ is some
standard numeration of $\EA$. We call such numerated theories \emph{G\"odelian extensions of $\EA$}.

The relation $\leq_\EA$ defines a natural preorder on the set $\fG_\EA$ of G\"odelian extensions of $\EA$. Let $\ol\fG_\EA$ denote the quotient by the associated equivalence relation $=_\EA$, where by definition $\gs=_\EA \tau$ iff both $\gs\leq_{\EA}\tau$ and $\tau\leq_{\EA}\gs$. $\ol\fG_\EA$ is a lattice with $\land_\EA$ corresponding to the union of theories and $\lor_\EA$ to their intersection. These operations are defined on elementary numerations as follows: \begin{eqnarray*}\gs\land_\EA\tau & := & \gs(x)\lor\tau(x),\\ \gs\lor_\EA\tau & := & \ex{x_1,x_2\leq x}(\gs(x_1)\land\tau(x_2)\land x=\text{disj}(x_1,x_2)),\end{eqnarray*}
where $\text{disj}(x_1,x_2)$ is an elementary term computing the G\"odel number of the disjunction of formulas given by G\"odel numbers $x_1$ and $x_2$.

We will only be concerned with the operation $\land_\EA$, that is, with the structure of lower semilattice with top $(\ol\fG_\EA,\land_\EA,1_\EA)$. Notice that the top element $1_\EA$ corresponds to (the equivalence class of) $\EA$, whereas the bottom $0_\EA$ is the class of all inconsistent extensions of $\EA$.

An operator $R:\fG_\EA\to \fG_\EA$ is called \emph{extensional} if $\gs=_\EA\tau$ implies $R(\gs)=_\EA R(\tau)$. Similarly, $R$ is called \emph{monotone} if $\gs\leq_\EA\tau$ implies $R(\gs)\leq_\EA R(\tau)$. Clearly, each monotone operator is extensional and each extensional operator correctly acts on the quotient lattice $\ol\fG_\EA$. An operator $R$ is called \emph{semi-idempotent} if $R(R(\gs))\leq_\EA R(\gs)$. $R$ is a \emph{closure} operator if it is monotone, semi-idempotent and, in addition, $\gs\leq_\EA R(\gs)$.  Operators considered in this paper will usually be at least monotone and semi-idempotent.

Meaningful monotone operators abound in arithmetic. Typical examples are the uniform $\Sigma_n$-reflection principles $\mR_n(\gs)$ associating with $\gs$ the extension of $\EA$ by the schema $\{\al{x}(\Box_\gs\phi(\bar x)\to \phi(x)):\phi\in\iPi_{n+1}\}$ taken with its natural elementary numeration that we denote $x\in \mR_n(\gs)$. It is known that the theory  $\mR_n(\gs)$ is finitely axiomatizable. Moreover, $\mR_0(\gs)$ is equivalent to G\"odel's consistency assertion $\Con(\gs)$ for $\gs$. The following basic lemma will be useful later.

Let $S$ be a G\"odelian extension of $\EA$ numerated by $\gs$, and let $x\in\Pi^0_n$ denote an elementary formula expressing that $x$ is the G\"odel number of a $\Pi^0_n$-sentence.
\bl \label{scompl} \benr
\item If $S$ extends $\EA$ by $\Pi^0_{n+1}$-axioms, then $\mR_n(\gs)$ contains $S$.
\item If $\EA\vdash \al{x}(\gs(x)\to x\in \iPi^0_{n+1})$ then $\mR_n(\gs)\leq_\EA \gs$.
\eenr
\el

\bp\ The second claim is a straightforward formalization of the first one. To prove Claim (i) assume $S\vdash \phi$. Then there is a $\pi\in\iPi^0_{n+1}$ such that $\EA\vdash \pi\to \phi$ and $S\vdash \pi$. We have $\EA\vdash \Box_\gs \pi$ by $\Sigma_1$-completeness. Then $\mR_n(\gs)\vdash\pi\vdash \phi$. \ep

In this paper we will study another series of monotone operators. Given a theory $S$ numerated by $\gs$, let $\iPi_n(S)$ denote the extension of $\EA$ by all theorems of $S$ of complexity $\Pi^0_n$. The set $\iPi_n(S)$ is r.e.\ but in general not elementary recursive. In order to comply with our definitions we apply a form of Craig's trick that yields an elementary axiomatization of $\iPi_n(S)$.\footnote{Over $\EA+\BS_1$ one can  work with a natural r.e.\ axiomatization of $\iPi_n(S)$.} Let $\iPi_n(\gs)$ denote the elementary formula $$\ex{y,p\leq x}(\Prf_\gs(y,p)\land y\in \Pi^0_n \land x=\text{disj}(y,\gn{\bar p\neq\bar p}))$$ and the theory numerated by this formula over $\EA$. Here, $\Prf_\gs(y,p)$ is an elementary formula expressing that $p$ is the G\"odel number of a proof of $y$, so that $\ex{p}\Prf_\gs(y,p)$ is $\Box_\gs(y)$. Then it is easy to see that the theory $\iPi_n(\gs)$ is (externally) deductively equivalent to $\iPi_n(S)$.

We will implicitly rely on the following characterization.

\bl \label{vis} It is provable in $\EA$ that
$$\al{x}(\Box_{\iPi_n(\gs)}(x)\eqv \ex{\pi\in\Pi^0_n}(\Box_\gs(\pi)\land\Box_\EA(\pi\to x)).$$
\el
\bp\ The implication from right to left is easy, we sketch a proof of $(\to)$. Reason within $\EA$. Suppose $p$ is a $\iPi_n(\gs)$-proof of $x$. It is a $\EA$-proof of $x$ from some assumptions $\pi'_1,\pi'_2,\dots,\pi'_k$ such that each $\pi'_i$ has the form $\pi_i\lor \ol{p_i}\neq \ol{p_i}$ where $\pi_i\in\Pi^0_n$ and $\Prf_\gs(\pi_i,p_i)$. Since $p$ contains witnesses for all the proofs $p_i$, from $p$ one can construct in an elementary way a sentence $\pi\in\Pi^0_n$ equivalent to $\pi_1\land\dots\land\pi_k$ together with its $\gs$-proof and an $\EA$-proof of $\pi\to x$, using a formalization of the deduction theorem in $\EA$. A verification that it is, indeed, the required proof goes by an elementary induction on the length of $p$.
\ep

Using Lemma \ref{vis} one can naturally infer that all the operators $\mR_n$ and $\iPi_n$ are monotone and semi-idempotent, moreover $\iPi_n$ is a closure. It is easy to see that $\EA$ can be replaced in all the previous considerations by any of its G\"odelian extensions $T$. The main source of interest for us in this paper will be the structure of semilattice with operators $$(\ol{\fG}_T,\land_T,1_T,\{\mR_n,\iPi_{n+1}:n<\gw\})$$ and its subsemilattice with operators $\ol{\fG}^0_T$ generated by $1_T$. We call the former \emph{the $\Rcn$-algebra of G\"odelian extensions of $T$}. The term \emph{$\Rcn$-algebra} will be explained below.

\section{Strictly positive logics and reflection calculi}

We refer the reader to a note \cite{Bek18b} for a short introduction to strictly positive logic sufficient for the present paper and to \cite{KKTWZ-arx} for more information from the general  algebraic perspective. For a general background on modal logic and provability logic we refer to the books~\cite{ChZa,Smo85,Boo93}.

\subsection{Normal strictly positive logics}
Consider a modal language $\cL_\gS$ with propositional variables
$p,q$,\dots , a constant $\top$, conjunction $\land$, and a possibly
infinite set of symbols $\Sigma=\{a_i:i\in J\}$ understood as
diamond modalities. The family $\Sigma$ is called the
\emph{signature} of the language $\cL_\gS$. Strictly positive
formulas (or simply \emph{formulas}) are built up by the grammar:
$$A::= p \mid \top \mid (A\land A) \mid a A, \quad \text{where $a\in \Sigma$.}$$
\emph{Sequents} are expressions of the form $A\vdash B$ where $A,B$
are strictly positive formulas.

Basic sequent-style system, denoted $\Kp$, is given by the following
axioms and rules:

\ben
\item $A\vdash A; \quad A\vdash\top; \quad$ if $A\vdash B$ and $B\vdash C$ then $A\vdash C$;
\item $A\land B\vdash A; \quad A\land B\vdash B; \quad$ if $A\vdash B$ and $A\vdash C$ then
$A\vdash B\land C$;
\item  if $A\vdash B$ then $a A\vdash a B$, for each $a\in\gS$.
\een

It is well-known that $\Kp$ axiomatizes the strictly positive fragment of a polymodal version of basic modal logic K (cf~\cite{Bek18b,KKTWZ-arx}). All our systems will also contain the following principle corresponding to the transitivity axiom in modal logic:

\ben \setcounter{enumi}{3}
\item $aa A\vdash a A$.
\een
The extension of $\Kp$ by this axiom will be denoted $\Kxp$~\cite{Das12}.

Let $C[A/p]$ denote the result of replacing in $C$ all occurrences
of a variable $p$ by $A$. A set of sequents $L$ is called a \emph{normal strictly positive logic} if it contains the axioms and is closed under the rules of $\Kp$ and under the
following \emph{substitution rule}: if $(A\vdash B)\in L$ then
$(A[C/p]\vdash B[C/p])\in L$. We will only consider normal strictly positive logics below. We write $A\vdash_L B$ for the statement that $A\vdash B$ is provable in $L$ (or belongs to $L$). $A=_L B$ means $A\vdash_L B$ and $B\vdash_L A$.

Any normal strictly positive logic $L$ satisfies the following simple \emph{positive replacement lemma} that we leave without proof.

\bl Suppose $A\vdash_L B$, then $C[A/p]\vdash_L C[B/p]$, for any formula
$C$.
\el

\subsection{Algebraic semantics}
Algebraic semantics for normal strictly positive logics is given by \emph{semilattices with monotone operators} (SLOs), that is, structures of the form $\fM=(M;\land^\fM,1^\fM,\{a^\fM: a\in \Sigma\})$ where $(M,\land^\fM,1^\fM)$ is a semilattice with top and each $a^\fM:M\to M$ is a monotone operator on $\fM$: $x\leq y$ implies $a^\fM(x)\leq a^\fM(y)$, for all $x,y\in M$. Every strictly positive formula $A$ of $\cL_\Sigma$ represents a term $A^\fM$ of $\fM$. We say that \emph{$A\vdash B$ holds in $\fM$} (or $\fM$ \emph{satisfies} $A\vdash B$) if $\fM\models\al{\vec x}A^\fM(\vec x)\leq B^\fM(\vec x)$. It is easy to see that $A\vdash_{\mathrm{K}^+} B$ if and only if $A\vdash B$ holds in each SLO $\fM$.
The SLOs satisfying all the theorems of a normal s.p.\ logic $L$ are called \emph{$L$-algebras}.

Given a normal s.p.\ logic $L$ in a signature $\Sigma$ and an alphabet of variables $V$, its \emph{Lindenbaum--Tarski algebra} is a SLO $\fL_L^V$ whose domain consists of the equivalence classes of formulas of the language of $L$ modulo $=_L$. Let $[A]_L$ denote the equivalence class of $A$. The operations are defined in a standard way as follows: $[A]_L\land^\fL [B]_L:= [A\land B]_L$, $a^\fL([A]_L):=[aA]_L$, for each $a\in\Sigma$. It is well-known that $A\vdash_L B$ iff $A\vdash B$ holds in $\fL_L^V$. Hence, any normal s.p.\ logic $L$ is complete w.r.t.\ its algebraic semantics, that is, w.r.t.\ the class of all $L$-algebras.

The algebra $\fL_L^V$ is also called the \emph{free $V$-generated $L$-algebra}. In this paper we will be particularly interested in the algebras $\fL_L^V$ where $V$ is empty. In this case we denote the algebra $\fL_L^V$ by $\fL_L^0$.


\subsection{The system $\Rc$}
\emph{Reflection calculus $\Rc$} is a normal strictly positive logic formulated in the signature $\{\Diamond_n:n\in\gw\}$. It is obtained by adjoining to the axioms and rules of $\Kxp$ (stated for each $\Diamond_n$) the following principles:
\ben \setcounter{enumi}{4}
\item $\Diamond_n A\vdash \Diamond_m A$, for all $n>m$;
\item $\Diamond_n A\land \Diamond_m B\vdash \Diamond_n(A\land \Diamond_m B)$, for all
$n>m$. \label{six}
\een
We notice that $\Rc$ proves the following \emph{polytransitivity} principles:
$$\Diamond_n\Diamond_m A\vdash \Diamond_m A, \quad \Diamond_m\Diamond_n A\vdash \Diamond_m A, \quad \text{for each $m\leq n$.}$$
Also, the converse of Axiom \ref{six} is provable in $\Rc$, so that in fact we have \beq \label{sixeq} \Diamond_n(A\land \Diamond_m B)=_\Rc \Diamond_n A\land \Diamond_m B.\eeq

The system $\Rc$ was introduced in an equational logic format by Dashkov \cite{Das12}, the present formulation is from \cite{Bek12a}. Dashkov showed that $\Rc$ axiomatizes the set of all sequents $A\vdash B$ such that the implication $A\to B$ is provable in the polymodal logic GLP. Moreover, unlike GLP itself, $\Rc$ is polytime decidable (whereas GLP is \textsc{PSpace}-complete~\cite{Shap08}) and enjoys the finite frame property (whereas GLP is Kripke incomplete).

We recall a correspondence between variable-free $\Rc$-formulas and ordinals \cite{Bek04}.
Let $\Fo$ denote the set of all variable-free $\Rc$-formulas, and let $\Fo_n$ denote its restriction to the signature $\{\Diamond_i:i\geq n\}$, so that $\Fo=\Fo_0$. For each $n\in\gw$ we define binary relations $<_n$ on $\Fo$ by
$$A<_n B \iffdef B\vdash_\Rc \Diamond_nA.$$
Obviously, $<_n$ is a transitive relation invariantly defined on the equivalence classes w.r.t.\ provable equivalence in $\Rc$ (denoted $=_\Rc$). Since $\Rc$ is polytime decidable, so are both $=_\Rc$ and all of $<_n$.

An $\Rc$-formula without variables and $\land$ is called a \emph{word} (or a \emph{worm} in some treatments). In fact, any such formula syntactically is a finite sequence of letters $\Diamond_i$
(followed by $\top$). If $A,B$ are words then $AB$ will denote $A[\top/B]$, that is, the word corresponding to the concatenation of these sequences. $A\circeq B$ denotes the graphical identity of  formulas (words).

The set of all words will be denoted $\Wo$, and $\Wo_n$ will denote its restriction to the signature $\{\Diamond_i:i\geq n\}$.  The following facts are from \cite{Bek04,Bek12a}:
\bi
\item Every $A\in \Fo_n$ is $\Rc$-equivalent to a word in $\Wo_n$;
\item $(\Wo_n/{=_\Rc},<_n)$ is isomorphic to $(\ge_0,<)$.
\ei

Here, $\ge_0$ is the first ordinal $\ga$ such that $\gw^\ga=\ga$. Thus, the set $\Wo_n/{=}_{\Rc}$ is well-ordered by the relation $<_n$. The isomorphism can be established by an onto and order preserving function
$o_n:\Wo_n\to \ge_0$ such that, for all $A,B\in \Wo_n$, $$A=_\Rc B \iff o_n(A)=o_n(B).$$ Then $o_n(A)$ is the order type of $\{B\in\Wo_n:B<_n A\}/{=_\Rc}$.

The function $o(A):=o_0(A)$ can be inductively calculated as follows: If $A\circeq \Diamond_0^k\top$ then $o(A)=k$. If $A\circeq A_1 \Diamond_0 A_2 \Diamond_0\cdots \Diamond_0 A_n$, where all
$A_i\in \Wo_1$ and not all of them are empty, then
\[o(A)= \gw^{o( A_n^-)}+\cdots +\gw^{o( A_1^-)}.\]
Here, $B^-$ is obtained from $B\in \Wo_1$ replacing every $\Diamond_{m+1}$ by $\Diamond_m$. For $n>0$ and $A\in\Wo_n$ we let $o_n(A)=o_{n-1}(A^-)$.

\subsection{The system $\Rcn$}

\bd
The signature of $\Rc^\nab$ consists of modalities $\Diamond_n$ and $\nab_n$, for each $n<\gw$. The system $\Rc^\nab$ is a normal strictly positive logic given by the following axioms and rules, for all $m,n<\gw$:

\ben
\item $\Rc$ for $\Diamond_n$; $\Rc$ for $\nab_n$;
\item $A\vdash \nab_n A$;
\item $\Diamond_n A\vdash \nab_n A$;
\item $\Diamond_m\nab_n A\vdash \Diamond_m A$ \ if $m\leq n$; \item $\nab_n \Diamond_m A\vdash \Diamond_m A$ \ if $m\leq n$.
\een
\ed

As a basic syntactic fact about $\Rcn$ we mention the following useful lemma. For brevity, we often write $=$ for $=_{\Rc^\nab}$ and $\vdash$ for $\vdash_\Rcn$.

\bl \label{iden} The following are theorems of $\Rc^\nab$, for all $m<n$:
\benr
\item $\Diamond_n(A\land\nab_m B) = \Diamond_n A \land \Diamond_m B$;
\item $\nab_n(A\land \Diamond_m B) = \nab_n A \land \Diamond_m B$.
\eenr
\el

\bp\ (i) Part $(\,\vdash\,)$ follows from $\Diamond_n\nab_m B\vdash \Diamond_m B$.
Part $(\,\dashv\,)$ follows from $\Diamond_n A \land \Diamond_m B\vdash \Diamond_n(A\land \Diamond_m B)\vdash \Diamond_n(A\land \nab_m B)$ using positive replacement.

(ii) Part $(\,\vdash\,)$ follows from $\nab_n\Diamond_m B\vdash \Diamond_m B$.
Part $(\,\dashv\,)$ follows from $\nab_n A \land \Diamond_m B\vdash \nab_n A \land \nab_m \Diamond_m B\vdash \nab_n(A\land \nab_m \Diamond_m B)\vdash \nab_n(A\land \Diamond_m B)$ using Axiom \ref{six} for $\nab$ modalities, the fact that $\Diamond_m B=\nab_m \Diamond_m B$ and positive replacement.
\ep

A formula $A$ is called \emph{ordered} if no modality with a smaller index $i$ (be it $\Diamond_i$ or $\nab_i$) occurs in $A$ within the scope of a modality with a larger index $j>i$.

\bl\label{ord}
Every formula $A$ of $\Rcn$ is equivalent to an ordered one.
\el

\bp\ Apply equation \refeq{sixeq} of $\Rc$ for $\Diamond$ and for $\nabla$ modalities, and the identities of Lemma \ref{iden} from left to right, until the rules are not applicable to any of the subformulas of $A$. \ep

The intended arithmetical interpretation of $\Rcn$ maps strictly positive
formulas to G\"odelian theories in $\fG_T$ in such a way that
$\top$ corresponds to $T$, $\land$
corresponds to the union of theories, $\Diamond_n$ corresponds to $\mR_n$ and $\nab_n$ corresponds to $\iPi_{n+1}$, for each $n\in\gw$.

\bd An \emph{arithmetical interpretation in $\fG_T$} is a map $*$ from strictly positive
modal formulas to $\fG_T$ satisfying the following conditions for all $n\in \gw$:

\bi
\item $\top^*=T$; \quad $(A\land B)^*=(A^*\land_T B^*)$;
\item $(\Diamond_n A)^*=\mR_n(A^*)$;
\quad $(\nab_n A)^*=\iPi_{n+1}(A^*)$. \ei \ed

The following result shows, as expected, that every theorem of $\Rcn$ represents an identity of the structure  $(\ol{\fG}_T,\land_T, 1_T,\{\mR_n,\iPi_{n+1}:n<\gw\})$.

\bt\ \label{soundness} For any formulas $A,B$ of $\Rc^\nab$, if $A\vdash_{\Rcn} B$ then $A^*\leq_T B^*$, for all arithmetical interpretations $*$ in $\fG_T$.
\et

\bp\ A proof of Theorem \ref{soundness} is routine by induction on the length of the derivation. For the axioms and rules of $\Rc$ for the $\Diamond$-fragment the claim has been carefully verified in \cite{Bek14}. The RC-axioms for the $\nab$-fragment are obvious except for Axiom \ref{six}, that is, the principle
\beq \nab_n A\land \nab_m B\vdash \nab_n(A\land \nab_m B).\label{sixn} \eeq
Consider any arithmetical interpretation $*$, and let $S=A^*$ and $U=B^*$ be the corresponding G\"odelian theories (with the associated numerations $\gs$ and $\tau$, respectively). We rely on Lemma~\ref{vis}. The principle \refeq{sixn} is the formalization in $\EA$ of the following assertion: \emph{For any sentence $\pi\in \Pi^0_{n+1}$, if $S\cup \iPi_{m+1}(U)\vdash \pi$ then $\iPi_{n+1}(S)\cup \iPi_{m+1}(U)\vdash \pi$.} Reasoning in $\EA$, consider a sentence $\phi\in \iPi_{m+1}(U)$ such that $S,\phi\vdash \pi$. Then $S\vdash \phi\to \pi$ and, since $\phi\to \pi$ is logically equivalent to a $\Pi^0_{n+1}$-sentence, conclude $\iPi_{n+1}(S)\vdash \phi\to \pi$. Thus, $\iPi_{n+1}(S)\cup \iPi_{m+1}(U)\vdash \pi$.

Concerning the remaining axioms of $\Rcn$ we remark that Axiom 2 holds since the theory $\iPi_{n+1}(S)$ is (provably) contained in $S$.
Axiom 3 is Lemma~\ref{scompl}~(ii).

Axiom 4: Assume $\mR_m(\iPi_{n+1}(\gs))$. In order to prove $\mR_{m}(\gs)$ let  $\phi\in\iPi_{m+1}$ and $\Box_\gs\phi$. Then clearly $\Box_{\iPi_{n+1}(\gs)}\phi$, since $m\leq n$, and hence $\phi$, by $\mR_m(\iPi_{n+1}(\gs))$.

Axiom 5 formalizes the fact that $\mR_{m}(\gs)$ is an extension of $T$ by a   $\Pi_{m+1}$-sentence.
\ep

Theorem \ref{soundness}, together with G\"odel's second incompleteness theorem, has as its corollary the following property of the logic $\Rcn$.
\bcor \label{irrefl} For all $\Rcn$ formulas $A$, $A\nvdash_{\Rcn} \Diamond_n A$.
\ecor

\bp\ Assume otherwise, then interpreting $\Rcn$ in $\fG_\EA$ yields $A^*\leq_\EA \mR_n(A^*)$ by Theorem \ref{soundness}. Hence, by G\"odel's theorem the theory $A^*$ is inconsistent. This contradicts the soundness of $\EA$. \ep

A similar fact is known for $\Glp$ and can also be proved by purely modal logic means \cite{BJV,BFJ14}. An elementary argument for $\Rc$ is given in Appendix A. David Fern\'andez-Duque gives another proof for a generalization of $\Rc$ with transfinitely many modalities. We will make use of Corollary~\ref{irrefl} (for $\Rc$) in the normal form theorems below. Whereas a reference to the given proof of Corollary~\ref{irrefl} presupposes at least the soundness of $\EA$, the elementary Kripke model argument for $\Rc$ is formalizable in $\EA$.

\begin{conjecture} $\Rcn$ is arithmetically complete, that is, the converse of Theorem~\ref{soundness} also holds, provided $T$ is arihmetically sound.
\end{conjecture}

\subsection{Kripke incompleteness of $\Rc^\nab$}

Kripke frames and models are understood in this paper in the usual
sense. A \emph{Kripke frame} $\cW$ for the language of $\Rc^\nab$ consists
of a non-empty set $W$ equipped with a family of binary relations
$\{R_n,S_n:n\in\gw\}$.

A \emph{Kripke model} is a Kripke frame $\cW$ together with a
\emph{valuation} $v: W\times \Var \to \{0,1\}$ assigning a truth
value to each propositional variable at every node of $\cW$. As
usual, we write $\cW,x\fc A$ to denote that a formula $A$ is true at
a node $x$ of a model $\cW$. This relation is inductively defined as
follows:

\bi
\item $\cW,x\fc p \iff v(x,p)=1$, for each $p\in\Var$;
\item $\cW,x\fc \top$; \quad $\cW,x\fc A\land B \iff (\cW,x\fc
A \text{ and } \cW,x\fc B)$;
\item $\cW,x\fc \Diamond_n A \iff \ex{y}(xR_n y\text{ and }
\cW,y\fc A)$;
\item $\cW,x\fc \nab_n A \iff \ex{y}(xS_n y\text{ and }
\cW,y\fc A)$.
\ei

A formula $A$ is \emph{valid in a Kripke frame $\cW$} if $\cW,x\fc A$, for each $x\in W$ and each valuation $v$ on $\cW$. The following lemma is standard and easy.

\bl
A Kripke frame $\cW$ validates all theorems of $\Rcn$ iff the following conditions hold, for all $m,n<\gw$:
\benr
\item $R_n$ is transitive; $R_n\subseteq R_m$ if $m<n$; $R_n^{-1}R_m\subseteq R_m$ if $m<n$;
\item $S_n$ is transitive, reflexive; $S_n\subseteq S_m$ if $m<n$; $S_n^{-1}S_m\subseteq S_m$ if $m<n$;
\item  $R_n\subseteq S_n$; $S_n R_m\subseteq R_m$, $R_mS_n\subseteq R_m$ if $m\leq n$.
\eenr
\el

By the following proposition $\Rcn$ turns out to be incomplete w.r.t.\ its Kripke frames.

\bpr
The sequent $$\Diamond_1A\land \nab_0 B\vdash \Diamond_1(A\land \nab_0 B) \eqno (*)$$ is valid in every Kripke frame satisfying $\Rcn$. However, it is unprovable in $\Rcn$ (and arithmetically invalid).
\epr
\bp\
Firstly, it is easy to see that conditions $R_1\subseteq S_1$ and $S_1^{-1}S_0\subseteq S_0$ imply $R_1^{-1}S_0\subseteq S_0$. Therefore, $(*)$ holds in each Kripke frame of $\Rcn$.

Secondly, take $\top$ for $A$ and $\Diamond_1\top$ for $B$. The left hand side is $\Rcn$-equivalent to $\Diamond_1\top$. The right hand side is equivalent to $\Diamond_1(\top\land \nab_0\Diamond_1\top)=_\Rcn\Diamond_1\top \land \Diamond_0\Diamond_1\top$, by Lemma \ref{iden}(i).
By Corollary \ref{irrefl},  $\Diamond_1\top\nvdash_{\Rcn}\Diamond_0\Diamond_1\top$. Hence, $(*)$ is unprovable in $\Rcn$. \ep

By Theorem 3 of \cite{Bek18b}, a normal strictly positive logic is a fragment of some normal modal logic if and only if it is Kripke complete. Hence, we obtain

\bcor $\Rcn$ is not a strictly positive fragment of any normal modal logic. \ecor

\section{The variable-free fragment of $\Rcn$} \label{var-free}

Let $\Fo_n^{\nab}$ denote the set of all variable-free strictly positive formulas in the language of $\Rcn$ with the modalities $\{\Diamond_i,\nab_i:i\geq n\}$ only, and let $\Fon$ denote $\Fo_0^{\nab}$.
We abbreviate $F\vdash_{\Rcn} \nab_n G$ by $F\vdash_n G$ and $\nab_n F =_\Rcn \nab_n G$ by $F\equiv_n G$. First, we are going to establish a crucial result that every formula in $\Fo_n^{\nab}$ is equivalent to a word in $\Wo_n$ modulo $\equiv_n$. From this fact we will infer a weak normal form theorem for the variable-free fragment of $\Rcn$. Second, we will obtain two different unique normal form theorems for the variable-free fragment by sharpening the weak normal forms.

\subsection{Weak normal forms}
We begin with a few auxiliary lemmas.
\bl \label{cons}
\benr\item
If $A\vdash_n B$ and $m<n$, then $A\land \Diamond_m C\vdash_n B\land \Diamond_m C$;
\item If $A\vdash_n B$ and $B\vdash \nab_n C$, then $A\vdash \nab_n C$;
\item If $A\vdash_n B$ and $B\vdash \Diamond_n C$, then $A\vdash \Diamond_n C$.
\eenr
\el

\bp\ (i) $A\land \Diamond_m C\vdash \nab_n B \land \Diamond_m C \vdash \nab_n (B\land \Diamond_m C)$.

(ii) $A\vdash \nab_n B\vdash \nab_n\nab_n C\vdash \nab_n C$;

(iii) $A\vdash \nab_n B\vdash \nab_n\Diamond_n C\vdash \Diamond_n C$. \ep

\bl \label{id-up}
\benr
\item $\Diamond_i A\land B = \nab_i(\Diamond_i A\land B) \land B$;
\item $\nab_i A\land B = \nab_i(\nab_i A\land B) \land B$.
\eenr
\el

\bp\
In both (i) and (ii) the implication
$(\vdash)$ follows from the axiom $C\vdash \nab_i C$. For $(\dashv)$ we obtain $\nab_i(\Diamond_i A\land B)\vdash \nab_i\Diamond_i A = \Diamond_i A$ for (i) and simlarly $\nab_i(\nab_i A\land B)\vdash \nab_i\nab_i A = \nab_i A$ for (ii). \ep

\bl \label{lino} The set of all formulas  $\{\Diamond_n F, \nab_n G: F,G\in \Wo_n\}$ is linearly ordered by $\vdash_\Rcn$. \el
\bp\ For any $F,G\in \Wo_n$ we know that either $F\vdash_\Rc \Diamond_n G$ or $G\vdash_\Rc \Diamond_n F$ or $F=_\Rc G$. In the first case we obtain provably in $\Rcn$: $\Diamond_n F\vdash\nab_n F\vdash \Diamond_n G\vdash \nab_n G$. The second case is symmetrical. In the third case we obtain $\Diamond_n F =\Diamond_n G \vdash \nab_n F=\nab_n G$. \ep

\bt \label{equiv} For each $A\in\Fo_n^{\nab}$ there is a word $W\in \Wo_n$ such that $A\equiv_n W$. \et

\bp\ By Lemma \ref{ord} it is sufficient to prove the theorem for ordered formulas $A$. The proof goes by induction on the length of ordered $A$. We can also assume that the minimal modality occurring in $A$ is $\Diamond_n$ or $\nab_n$. (Otherwise, prove it for the minimum $m>n$ and infer $A\equiv_n W$ from $A\equiv_m W$.) The basis of induction is trivial, we consider the induction step.

Assume that the induction hypothesis holds for all formulas shorter than $A$. Since $A$ is ordered, $A$ can be written in the form $$A=\Diamond_n A_1\land \dots \land \Diamond_n A_k\land \nab_n B_1\land \dots \nab_n B_l\land D,$$
where $D\in \Fo_{n+1}^\nab$ and $A_i,B_j\in \Fo_n^{\nab}$. Since $\Diamond_n$ or $\nab_n$ must occur in $A$, we know that $D$ and each $A_i,B_j$ are strictly shorter than $A$. By the induction hypothesis and Lemma \ref{lino} we can delete from the conjunction all but one members of the form $\Diamond_n A_i$, $\nab_n B_j$. Thus, $A= D\land \Diamond_n A'$ or $A= D\land \nab_n B'$, for some words $A', B'\in \Wo_n$.

Now we apply the induction hypothesis to $D$ and obtain a word $V\in \Wo_{n+1}$ such that $V\equiv_{n+1} D$. It follows that $D\land \Diamond_n A'\equiv_{n+1} V\land \Diamond_n A'$ and $D\land \nab_n B'\equiv_{n+1} V\land \nab_n B'$, by Lemma \ref{cons}. Hence, it is sufficient to prove that, for some $W\in \Wo_n$, $V\land \Diamond_n A'\equiv_n W$ and similarly, for some $W\in \Wo_n$, $V\land \nab_n B'\equiv_n W$.

In the first case we actually have $V\land \Diamond_n A' =_\Rc W$, for some word $W$, which immediately yields the claim.

In the second case we write $B'=B_1\Diamond_n B_2$ where $B_1\in \Wo_{n+1}$. There are three cases to consider: (a) $B_1\vdash\Diamond_{n+1} V$, (b) $V\vdash \Diamond_{n+1} B_1$, (c) $V=B_1$.

In case (c) by Lemma \ref{id-up} we obtain:
$$V\land \nab_n B_1\Diamond_n B_2 = V\land \nab_n(V\land \Diamond_n B_2) = V\land \Diamond_n B_2 = V\Diamond_n B_2.$$

In case (a) we show $\nab_n(V\land \nab_n B')=\nab_n B'$. Firstly, $$B'\vdash \Diamond_{n+1}V\land \nab_n B'\vdash \nab_{n+1} V\land \nab_n B'=\nab_{n+1}(V\land \nab_n B').$$ Hence, $\nab_n B'\vdash \nab_n\nab_{n+1}(V\land \nab_n B')=\nab_n(V\land \nab_n B').$ On the other hand, $$\nab_n(V\land \nab_n B')\vdash \nab_n\nab_n B'\vdash \nab_n B'.$$

In case (b) we show $\nab_n(V\land \nab_n B')=\nab_n(V\land \Diamond_n B_2)$ so that one can infer $\nab_n(V\land \nab_n B')=\nab_n V\Diamond_n B_2$.
On the one hand, we have $$\nab_n B'=\nab_n (B_1\land \Diamond_n B_2)\vdash \nab_n\Diamond_n B_2=\Diamond_n B_2,$$ which implies $\nab_n(V\land \nab_n B')\vdash \nab_n(V\land \Diamond_n B_2)$. On the other hand,
$$V\land \Diamond_n B_2= V\land \Diamond_{n+1}B_1 \land \Diamond_n B_2= V\land \Diamond_{n+1}(B_1\land \Diamond_n B_2)=V\land \Diamond_{n+1} B'\vdash V\land \nab_n B'.$$ Hence,  $\nab_n(V\land \Diamond_n B_2)\vdash\nab_n(V\land \nab_n B')$. \ep

From Theorem \ref{equiv} we obtain the following strengthening of Lemma \ref{lino}.
\bcor
\label{lino1} The set of all formulas  $\{\Diamond_n F, \nab_n G: F,G\in \Fo_n^{\nab}\}$ is linearly ordered by $\vdash_\Rcn$.
\ecor

\bcor
\label{lino2} For all formulas $A,B\in\Fo_n^{\nab}$, either $A\vdash \Diamond_n B$, or $B\vdash \Diamond_n A$, or $A\equiv_n B$.
\ecor
\bp\ Consider the words $A_1\equiv_n A$ and $B_1\equiv_n B$. By the linearity property for words either $A_1\vdash \Diamond_n B_1$ or $B_1\vdash \Diamond_n A_1$ or $A_1=B_1$. In the first case we obtain $A\vdash \nab_n A_1\vdash \nab_n\Diamond_n B_1\vdash \Diamond_n B_1\vdash \Diamond_n \nab_n B\vdash \Diamond_n B$. The second case is symmetrical, the third one implies $A\equiv_n B$ immediately. \ep

\bcor
\label{lino3} For all $A,B\in\Fo_n^{\nab}$, $\Diamond_n A\vdash \Diamond_n B$ iff $A\vdash \nab_n B$.
\ecor
\bp\ Assume $\Diamond_n A\vdash \Diamond_n B$. By Corollary \ref{lino2}, either $A\vdash \Diamond_n B$, or $B\vdash \Diamond_n A$, or $A\equiv_n B$. In the first and the third cases we immediately have $A\vdash \nab_n B$. In the second case  we obtain $\Diamond_n A\vdash \Diamond_n B\vdash \Diamond_n\Diamond_n A$ contradicting Corollary \ref{irrefl}.

In the opposite direction, if $A\vdash \nab_n B$ then $\Diamond_n A\vdash \Diamond_n\nab_n B\vdash \Diamond_n B$. \ep

\bt[weak normal forms] Every formula $A\in\Fo_n^{\nab}$ is equivalent in $\Rcn$ to a formula of the form
$$\nab_nA_n\land\nab_{n+1}A_{n+1}\land\dots\land\nab_{n+k}A_{n+k},$$ for some $k$, where $A_i\in \Wo_i$ for all $i=n,\dots, n+k$.
\et

\bp\ Induction on the build-up of $A\in\Fo_n^{\nab}$. We consider the following cases.

1) $A=B\land C$. The induction hypothesis is applicable to $B$ and $C$, so it is sufficient to prove: for any $B_i,C_i\in \Wo_i$ there is a word $A_i\in \Wo_i$ such that
$$\nab_i B_i\land \nab_i C_i=\nab_i A_i.$$
By Lemma \ref{lino} we can take one of $B_i,C_i$ as $A_i$.

2) $A=\nab_i B$, for some $i\geq n$. Then we obtain
\begin{multline*}
\nab_i B = \nab_i(\nab_nB_n\land\nab_{n+1}B_{n+1}\land\dots\land\nab_{n+k}B_{n+k}) =  \\
= \nab_nB_n\land\dots \land \nab_{i-1}B_{i-1}\land \nab_i(\nab_i B_i\land\dots\land\nab_{n+k}B_{n+k}) =\\
= \nab_nB_n\land\dots \land \nab_{i-1}B_{i-1}\land \nab_i B'_i,
\end{multline*}
for some $B'_i\in \Wo_i$, by Theorem \ref{equiv}.

3) $A=\Diamond_i B$, for some $i\geq n$. Then we obtain, using Lemma \ref{iden},
\begin{multline*}
\Diamond_i B = \Diamond_i(\nab_nB_n\land\nab_{n+1}B_{n+1}\land\dots\land\nab_{n+k}B_{n+k}) =  \\
=\Diamond_nB_n\land\dots \land \Diamond_{i-1}B_{i-1}\land \Diamond_i(\nab_i B_i\land\dots\land\nab_{n+k}B_{n+k}) =\\
=\nab_n\Diamond_nB_n\land\dots \land \nab_{i-1}\Diamond_{i-1}B_{i-1}\land\nab_i\Diamond_i B'_i,
\end{multline*}
for some $B'_i\in \Wo_i$, by Theorem \ref{equiv}. \ep

Weak normal forms are, in general, not unique. However, the following lemma and its corollary show that the ``tails'' of the weak normal forms are invariant (up to equivalence in $\Rcn$).

\bl \label{invar}
Let $A\circeq\nab_nA_n\land\nab_{n+1}A_{n+1}\land\dots\land\nab_{k}A_{k}$  and $B\circeq\nab_nB_n\land\nab_{n+1}B_{n+1}\land\dots\land\nab_{m}B_{m}$  be weak normal forms, $B_m\not\circeq \top$ and $A\vdash B$. Then $k\geq m$ and for all $i$ such that $n\leq i\leq k$ there holds
\benr
\item $\nab_iA_i\land\dots\land\nab_{k}A_{k}\vdash_{i} \nab_iB_i\land\dots\land\nab_{m}B_{m};$
\item $\nab_iA_i\land\dots\land\nab_{k}A_{k}\vdash \nab_iB_i\land\dots\land\nab_{m}B_{m}.$
\eenr
\el

\bp\ By definition, Claim (ii) implies Claim (i), but we first prove (i) and then strengthen it to (ii). For $i=n$ both claims are vacuous, so we assume $i>n$.

Denote $\ol{A_i}:=\nab_iA_i\land\dots\land\nab_{k}A_{k}$ and
$\ol{B_i}:=\nab_iB_i\land\dots\land\nab_{m}B_{m}$. By Lemma \ref{lino2} we have either $\ol{A_i}\vdash \Diamond_i \ol{B_i}$ or $\ol{B_i}\vdash \Diamond_i \ol{A_i}$ or $\ol{A_i}\equiv_i\ol{B_i}$. In the first and in the third case we obviously have $\ol{A_i}\vdash_i\ol{B_i}$ as required.

Assume $\ol{B_i}\vdash \Diamond_i \ol{A_i}$. Consider the formula $$C:= \Diamond_n A_n\land \dots \land \Diamond_{i-1} A_{i-1}\land \ol{B_i}.$$ We show that $C\vdash \Diamond_i C$ contradicting Corollary \ref{irrefl}.

Using our assumption and Lemma \ref{iden} (i) we obtain
\begin{eqnarray*}
C & \vdash & \Diamond_n A_n\land \dots \land \Diamond_{i-1} A_{i-1}\land \Diamond_i \ol{A_i} \\
& \vdash & \Diamond_i (\nab_n A_n\land \dots \land \nab_{i-1} A_{i-1}\land \ol{A_i}) \\
& \vdash &  \Diamond_n A_n\land \dots \land \Diamond_{i-1} A_{i-1}\land \Diamond_i A \\
& \vdash &  \Diamond_n A_n\land \dots \land \Diamond_{i-1} A_{i-1}\land \Diamond_i B \\
& \vdash & \Diamond_i ( \Diamond_n A_n\land \dots \land \Diamond_{i-1} A_{i-1}\land B) \\
& \vdash & \Diamond_i C.
\end{eqnarray*}
This proves Claim (i).

To prove (ii) assume the contrary and consider the maximal number $i$ such that $\ol{A_i}\nvdash \ol{B_i}$. Such an $i$ exists, since both $A$ and $B$ have finitely many terms. Thus, we have $\ol{A_{i+1}}\vdash \ol{B_{i+1}}$ and $$\nab_i A_i\land \ol{A_{i+1}}\nvdash \nab_i B_i\land \ol{B_{i+1}}.$$ It follows that $\nab_i A_i\land \ol{A_{i+1}}\nvdash \nab_i B_i=\nab_i\nab_i B_i$, hence $\ol{A_i}\nvdash_i \nab_i B_i$. Since $\ol{B_i}\vdash \nab_iB_i$, we obtain $\ol{A_i}\nvdash_i \ol{B_i}$ contradicting Claim (i).
\ep

\bcor \label{uni} Let $\nab_nA_n\land\nab_{n+1}A_{n+1}\land\dots\land\nab_{k}A_{k}$ be any  weak normal form of a formula $A\in\Fo_n^{\nab}$ with $A_k\not\circeq\top$. Then $k$ and each tail $\nab_iA_i\land\dots\land\nab_{k}A_{k}$ is defined uniquely up to equivalence in $\Rcn$.
\ecor

There are two formats for graphically unique normal forms. We call them `fat' and `thin', because the former consist of larger expressions, whereas the latter are obtained by pruning certain parts of a given formula. Fat normal forms, presented below, have a natural proof-theoretic meaning and are tightly related to collections of proof-theoretic ordinals called \emph{conservativity spectra} or \emph{Turing--Taylor expansions}~\cite{Joo15a}.

\subsection{Fat normal forms}

\bd \label{fnf} A formula $A\in\Fo_n^{\nab}$ is in the \emph{fat normal form for $\Fo_n^{\nab}$} if either $A\circeq\top$ or it has the form  $\nab_nA_n\land\nab_{n+1}A_{n+1}\land\dots\land\nab_{n+k}A_{n+k}$, where for all $i=n,\dots, n+k$, $A_i\in \Wo_i$, $A_{n+k}\not\circeq \top$ and
$$\nab_i A_i\vdash \nab_i(\nab_iA_i\land\dots\land\nab_{n+k}A_{n+k}).\eqno (*)$$
A variable-free formula $A$ is in the \emph{fat normal form} if $A$ is in the fat normal form for $\Fo_0^{\nab}$.
\ed

\brem \label{fnfo} In a fat normal form, for each $i$ such that $n\leq i\leq n+k$, there holds $\nab_i A_i=_{\Rcn} \nab_i(\nab_iA_i\land\dots\land\nab_{n+k}A_{n+k}).$
\erem

\bt \label{upnorm}
\benr \item Every $A\in\Fo_n^{\nab}$ is equivalent to a formula in the fat normal form for $\Fo_n^{\nab}$.
\item For any $A\in\Fo_n^{\nab}$, the words $A_i$ in the fat normal form of $A$ for $\Fo_n^{\nab}$ are unique modulo equivalence in $\Rc$.
\eenr
\et

\bp\ (i) First, we apply Theorem \ref{equiv}. Then, by induction on $k$ we show that any formula $\nab_nA_n\land\dots\land\nab_{n+k}A_{n+k}$ can be transformed into one satisfying $(*)$.

For $k=0$ the claim is trivial. Otherwise, by the induction hypothesis we can assume that $(*)$ holds for
$i=n+1,\dots, n+k$. Then we argue using Lemma \ref{id-up} as follows:
\begin{multline*}
\nab_nA_n\land\nab_{n+1}A_{n+1}\land\dots\land\nab_{n+k}A_{n+k}=  \\
=\nab_n(\nab_nA_n\land\nab_{n+1}A_{n+1}\land\dots\land\nab_{n+k}A_{n+k})\land \nab_{n+1}A_{n+1}\land\dots\land\nab_{n+k}A_{n+k} = \\
=\nab_n A'_n \land \nab_{n+1}A_{n+1}\land\dots\land\nab_{n+k}A_{n+k},
\end{multline*}
where $A'\in \Wo_n$ is obtained from Theorem \ref{equiv}.
Notice that
\begin{multline*}
\nab_n A'_n\vdash \nab_n(\nab_nA_n\land\nab_{n+1}A_{n+1}\land\dots\land\nab_{n+k}A_{n+k}) \vdash \\ \vdash
\nab_n(\nab_nA'_n\land\nab_{n+1}A_{n+1}\land\dots\land\nab_{n+k}A_{n+k}),
\end{multline*} hence $(*)$ holds for $i=n$. This proves Claim (i).

To prove Claim (ii) we apply Lemma \ref{invar}. Assume $A\vdash B$,
$A=\nab_nA_n\land\dots\land\nab_{n+k}A_{n+k}$ is in the fat normal form and $B=\nab_nB_n\land\dots\land\nab_{n+m}B_{n+m}$ is in a weak normal form. Then $k\geq m$ and, for all $i=n,\dots, n+m$, $\nab_{i} A_i\vdash \nab_{i} B_i$.

It follows that, if $A,B\in \Fo_n^{\nab}$ are both in the fat normal form and $A= B$ in $\Rcn$, then $m=k$ and $\nab_i A_i=\nab_i B_i$, for $i=n,\dots,n+k$. Since $\Wo_i$ is linearly preordered by $<_i$ in $\Rc$, the latter is only possible if $A_i=_\Rc B_i$. \ep

\brem As stated in Theorem \ref{upnorm}, fat normal forms are only unique modulo equivalence of the constituent words $A_i$ in $\Rc$. However, we know that words have graphically unique $\Rc$-normal forms~\cite{Bek04}. Combining the two notions together yields graphically unique normal forms for $\Rcn$.
\erem

Thus, we can test the equality of two variable-free formulas in $\Rcn$ by graphically comparing their unique normal forms. Alternatively, we observe the following property.

\bl  Let $A\circeq\nab_nA_n\land\nab_{n+1}A_{n+1}\land\dots\land\nab_{k}A_{k}$  and $B\circeq\nab_nB_n\land\nab_{n+1}B_{n+1}\land\dots\land\nab_{m}B_{m}$ be any fat normal forms. Then $A\vdash_\Rcn B$ holds iff $k\geq m$ and, for all $i$ such that $n\leq i\leq m$, one has $\Diamond_iA_i\vdash_\Rc \Diamond_i B_i$.
\el

\bp\ By Lemma~\ref{invar} and Remark~\ref{fnfo}, $A\vdash_\Rcn B$ holds iff, for all $i$ such that $n\leq i\leq m$, one has $\nab_iA_i\vdash_\Rcn \nab_iB_i$. However, the latter is equivalent to $A_i\vdash_\Rcn \nab_iB_i$ and to $\Diamond_i A_i\vdash_\Rcn \Diamond_i B_i$ by Corollary \ref{lino3}. Since words are linearly preordered in $\Rc$, the latter is also equivalent to $\Diamond_i A_i\vdash_\Rc \Diamond_i B_i$. \ep

The transformation of a variable-free formula to its fat normal form is computable. Hence, we obtain
\bcor The set of variable-free sequents $A\vdash B$ provable in $\Rcn$ is decidable.
\ecor

From the uniqueness of normal forms we also obtain arithmetical completeness of the variable-free fragment of $\Rcn$ in the standard way.

\bcor Suppose $A,B$ are variable-free and $T$ is a sound G\"odelian extension of $\EA$. Then $A\vdash_\Rcn B$ iff $A^*\leq_T B^*$, for all arithmetical interpretations $*$ in $\fG_T$.
\ecor

\bcor Suppose $T$ is a sound G\"odelian extension of $\EA$. Then the algebra $\ofG_T^0$ is isomorphic to the Lindenbaum--Tarski algebra of the variable-free fragment of $\Rcn$.
\ecor

\subsection{Thin normal forms} \label{thin-nf}
Let $A\circeq\nab_0A_0\land\nab_{1}A_{1}\land\dots\land\nab_{k}A_{k}$ be in a weak normal form. As before, we denote $\ol{A_i}\circeq\nab_{i}A_{i}\land\dots\land\nab_{k}A_{k}$.

\bd $A$ is in a \emph{thin normal form} if either $A\circeq\top$ or  $A\circeq\nab_0A_0\land\nab_{1}A_{1}\land\dots\land\nab_{k}A_{k}$ where $A_k\not\circeq \top$, for all $i< k$ $A_i\in \Wo_i$, and there is no $B_i\in \Wo_i$ such that $B_i<_i A_i$ and $\ol{A_i}=_{\Rcn}\nab_iB_i\land \ol{A_{i+1}}$.
\ed

This definition allows one to easily prove the existence and uniqueness of normal forms using the fact that words in $\Wo_i$ are pre-wellordered by $<_i$.

\bt\label{thin}
For each $A\in \Fon$ there is a thin normal form equivalent to $A$ in $\Rcn$. The thin normal form is unique modulo equivalence of the constitutent words $A_i$ in $\Rc$.
\et

\bp\ We recursively define the words $A_k, A_{k-1},\dots, A_0$. To determine $k$ and $A_k$ one takes any weak normal form for $A$ (observe that $\top$ is the $<_i$-minimum for each $i$). Once one has defined $A_k,\dots,A_{i+1}$ one can define $A_i$ by considering all the weak normal forms with the given $\nab_{i+1} A_{i+1}$, \dots, $\nab_k A_k$ and selecting the one with the $<_i$-minimal $A_i$.
By induction on $k-i$ it is also easy to see that all the words $A_k,\dots, A_0$ are thus uniquely determined modulo $\Rc$.  \ep

The given proof, though short, is non-constructive. Now we will show that the thin normal form can be effectively computed. First, we consider a particular case when the given weak normal form is $\nab_0 A \land \nab_1 B$. Then we will reduce the general case to this one.

Let $A\in \Wo_0, B\in \Wo_1, B\not\circeq\top$ and $A\circeq A_0\Diamond_0 A_1\Diamond_0\dots \Diamond_0 A_n$ with $A_i\in \Wo_1$.  If $B\vdash\nab_0 A$ then $\nab_0 A \land \nab_1 B=\nab_0 \top\land \nab_1 B$ which is its thin normal form. So we assume $B\nvdash \nab_0 A$.
We define
$$B| A := A_i\Diamond_0\dots \Diamond_0 A_n,$$ where $i$ is the least such that $B\leq_1 A_i$. Such an $i$ exists, for otherwise $B\vdash \Diamond_0A\vdash \nab_0 A$. Clearly, $B| A$ can be found effectively from $A$ and $B$ by deleting the appropriate initial segment of $A$. Also notice that $B\land A =_\Rc B\land (B| A)$. We consider three cases.

\textsc{Case 1:} $A_0>_1 B$. We claim that $\nab_0 A \land \nab_1 B$ is in a thin normal form.

Assume $A'<_0 A$, then $A\vdash \Diamond_1 B\land \Diamond_0 A'=\Diamond_1 (B\land \Diamond_0 A')\vdash \Diamond_1(\nab_0 A' \land \nab_1 B)$. Hence, if $\nab_0 A'\land \nab_1 B\vdash \nab_0 A$, then $A\vdash \Diamond_1(\nab_0 A \land \nab_1 B)\vdash \Diamond_1 A$ contradicting Corollary \ref{irrefl} .

\textsc{Case 2:} $A_0<_1 B$. We claim that $\nab_0 \Diamond_0 (B| A) \land \nab_1 B$ is the thin normal form of $\nab_0 A\land \nab_1 B$. Firstly, we show that $\Diamond_0 (B| A)\land \nab_1 B\vdash \nab_0 A$. By downwards induction on $j:=i$ to $0$ we show that
$$\Diamond_0 (B| A)\land \nab_1 B\vdash \Diamond_0 (A_j\Diamond_0\dots \Diamond_0A_n).$$
Basis of induction holds since $B| A=A_i\Diamond_0\dots \Diamond_0A_n$. Assume the claim holds for $j$. Since $\nab_1 B\vdash \nab_1 \Diamond_1 A_{j-1}=\Diamond_1 A_{j-1}$, we obtain:
\begin{eqnarray*}
\Diamond_0 (B|A)\land \nab_1 B & \vdash &  \Diamond_0 (A_j\Diamond_0\dots \Diamond_0A_n) \land \nab_1 B \\
& \vdash & \Diamond_1 A_{j-1}\land \Diamond_0 (A_j\Diamond_0\dots \Diamond_0A_n) \\
& \vdash & \Diamond_1 A_{j-1}\Diamond_0 A_j\Diamond_0\dots \Diamond_0A_n, \text{since $A_{j-1}\in S_1$} \\
& \vdash & \Diamond_0 A_{j-1}\Diamond_0 A_j\Diamond_0\dots \Diamond_0A_n.
\end{eqnarray*}
Hence, the claim holds for $j-1$ and by induction we conclude that $$\Diamond_0 (B|A)\land \nab_1 B\vdash \Diamond_0 A\vdash \nab_0 A.$$

Now we need to show that for all $A'<_0 \Diamond_0(B|A)$ one has $\nab_0 A'\land \nab_1 B\nvdash \nab_0 A$. If $\Diamond_0(B|A)\vdash \Diamond_0 A'$ then by Lemma \ref{lino3} $B|A\vdash \nab_0 A'$. Also $B|A\vdash A_i\vdash \nab_1 B$, since we assume $B\leq_1 A_i$. It follows that $B|A\vdash \nab_0 A'\land \nab_1 B$. On the other hand, $A\vdash \Diamond_0 (B|A)$ and $\nab_0A\vdash\Diamond_0 (B|A)$ whence $\nab_0 A'\land \nab_1 B\nvdash\nab_0A$ by Corollary \ref{irrefl}.

\textsc{Case 3:} $A_0=B$. Let $C:=A_1\Diamond_0\dots \Diamond_0 A_n$, thus $A\circeq B\Diamond_0 C$.We claim that $\nab_0\Diamond_0 C\land \nab_1B$ is the thin normal form of $\nab_0 A\land \nab_1 B$.

First, $\nab_1 B\land \Diamond_0 C\vdash \nab_1(B\land \Diamond_0 C)\vdash \nab_0(B\Diamond_0 C)=\nab_0 A.$ Hence, $\nab_0\Diamond_0 C\land \nab_1B=\nab_0 A\land \nab_1 B$.

Second, we show that if $A'<_0 \Diamond_0 C$ then $\nab_0 A'\land\nab_1 B\nvdash \nab_0 A$. Assume $A'<_0 \Diamond_0 C$. By Lemma \ref{lino3} we have $C\vdash \nab_0 A'$. Also, since $A_1\geq_1 B$, we have $C\vdash A_1\vdash \nab_1 B$ by Lemma \ref{lino3}. It follows that $C\vdash \nab_0 A'\land \nab_1 B$. On the other hand, $A=B\land \Diamond_0 C\vdash \Diamond_0 C$, hence $\nab_0 A\vdash \nab_0\Diamond_0 C\vdash \Diamond_0 C$. Therefore, by Corollary \ref{irrefl} $\nab_0 A'\land\nab_1 B\nvdash \nab_0 A$.

In all three cases we have explicitly constructed the thin normal form. Hence, we obtain the following theorem.

\bt
For any variable-free formula of $\Rcn$, its unique thin normal form can be effectively constructed.
\et

\bp\ Let a formula $A\circeq\nab_0A_0\land\nab_{1}A_{1}\land\dots\land\nab_{k}A_{k}$ in a weak normal form be given. We argue by induction on $k$. For $k=0$ the claim is obvious. Consider $k>0$, by IH we may assume that $\ol{A_1}:=\nab_{1}A_{1}\land\dots\land\nab_{k}A_{k}$ is in a thin normal form. (To formally apply the IH one should consider the formula obtained from $\ol{A_1}$ by decreasing all indices of modalities by $1$.)
By Theorem \ref{equiv} there is a word $B\in S_1$ such that $\nab_1 B\equiv_1\nab_{1}A_{1}\land\dots\land\nab_{k}A_{k}$.

Consider the formula $\nab_0A_0\land \nab_1 B$ and bring it to a thin normal form, that is, find a $<_0$-minimal $A_0'\in S_0$ such that $\nab_0A_0'\land \nab_1 B = \nab_0A_0\land \nab_1 B$. We claim that $A':=\nab_0A_0'\land\nab_{1}A_{1}\land\dots\land\nab_{k}A_{k}$ is equivalent to $A$ and is in a thin normal form.

Firstly, $A'\vdash \nab_0A_0'\land \nab_1 B\vdash \nab_0 A_0$, hence $A'\vdash A$. On the other hand, $A\vdash \nab_0A_0\land \nab_1 B = \nab_0A_0'$, hence $A\vdash A'$.

Secondly, assume there is an $A''<_0 A'_0$ such that $\nab_0 A'' \land \ol{A_1}\vdash \nab_0 A'_0$. By Lemma~\ref{cons}, $\nab_0 A'' \land \ol{A_1}\equiv_1 \nab_0 A'' \land \nab_1 B$. Hence,  $\nab_0 A'' \land \nab_1 B\vdash \nab_0 A'_0$ contradicting the $<_0$-minimality of $A'_0$. \ep

\section{Iterating monotone operators on $\fG_\EA$}

Transfinite iterations of reflection principles play an important role in proof theory starting from the works of A. Turing~\cite{Tur39} and S. Feferman~\cite{Fef62} on recursive progressions. Here we present a general result on defining iterations of monotone semi-idempotent operators in $\fG_\EA$.

An operator $R:\fG_\EA\to \fG_\EA$  is called \emph{computable} if so is the function $\gn{\gs}\mapsto \gn{R(\gs)}$. By extension of terminology we also call computable any operator $R'$ such that $\al{\gs\in \fG_\EA} R'(\gs)=_\EA R(\gs)$, for some computable $R$.

Bounded formulas in the language of $\EA$ will henceforth be called \emph{elementary}. An operator $R:\fG_\EA\to \fG_\EA$ is called \emph{uniformly definable} if there is an elementary formula $\Ax_R(x,y)$ such that
\benr \item For each $\gs\in \fG_\EA$ one has $R(\gs)=_\EA \Ax_R(x,\ol{\gn{\gs}})$,
\item $\EA\vdash \al{x,y}(\Ax_R(x,y) \to x\geq y).$
\eenr

The operators $\mR_n$ and $\iPi_{n+1}$ are uniformly definable in a very special way. For example, the formula $\mR_n(\gs)$ is obtained by substituting $\gs(x)$ for $X(x)$ into a fixed elementary formula containing a single positive occurrence of a predicate variable $X$. More generally, it can be shown that an operator $R:\fG_\EA\to \fG_\EA$ is uniformly definable iff $R$ is computable. A proof of this fact is given in Appendix B.

\bd
A uniformly definable $R$ is called
\bi
\item \emph{provably monotone} if $\EA\vdash \al{\gs,\tau}(\bb\tau\leq_\EA \gs\qq \to \bb R(\tau)\leq_\EA R(\gs)\qq),$
\item \emph{reflexively monotone} if $\EA\vdash \al{\gs,\tau}(\bb\tau\leq_\EA \gs\qq \to \bb R(\tau)\leq R(\gs)\qq).$
\ei  \ed
Here, $\gs,\tau$ range over G\"odel numbers of elementary formulas in one free variable, $\bb \tau\leq_\EA\gs\qq $ abbreviates $\Box_\EA\al{x}(\Box_\gs(x)\to \Box_\tau(x))$ and
$\bb R(\tau)\leq R(\gs)\qq $ stands for $\al{x}(\Box_{\Ax_R(\cdot,\bar\gs)}(x)\to \Box_{\Ax_R(\cdot,\bar\tau)}(x)).$ Reflexive monotonicity here refers to the fact that $\bb R(\tau)\leq R(\gs)\qq $ is the statement of inclusion of theories rather than provable inclusion. Since the formula $\bb \tau\leq_\EA\gs\qq $ implies its own provability in $\EA$, reflexively monotone operators are (provably) monotone but not necessarily vice versa.

It is also easy to see that the operators $\mR_n$ (along with all the usual reflection principles) are reflexively monotone.

\renewcommand{\preceq}{\preccurlyeq}

Next we turn to iterations of operators along ordinal notation systems. In this paper, ordinal notation systems will be \emph{pre-wellorderings}, that is, reflexive, transitive binary relations whose quotient order is a well-ordering.
An \emph{elementary pre-wellordering} is a pair of bounded formulas $D(x)$ and $x\preceq y$ and a constant $0$ such that the relation $\preceq$ provably in $\EA$ is a linear preorder on $D$ with the least element $0$, and is a pre-wellorder of $D$ in the standard model of arithmetic. Given an elementary well-ordering $(D,\preceq,0)$, we will denote its elements by Greek letters and will identify them with an initial segment of the ordinals. We denote \begin{eqnarray*}
x\approx y & \iffdef & (x\preceq y \land y\preceq x); \\
x\prec y & \iffdef & (x\preceq y \land y\not\preceq x).
\end{eqnarray*}

Let $R$ be an uniformly definable monotone operator. The \emph{$\ga$-th iterate of $R$ along $(D,\preceq,0)$} is a map associating with any numeration $\gs$ the  G\"odelian extension of $\EA$ numerated by an elementary formula $\rho(\ol\ga,x)$ such that provably in $\EA$:
\beq\rho(\ga,x)\eqv ((\ga\approx 0\land \gs(x))\lor \ex{\gb\prec\ga} \Ax_R(x,\gn{\rho(\bar \gb,x)})). \label{itref} \eeq
A natural G\"odel numbering of formulas and terms should satisfy the inequalities $\gn{\rho(\bar\gb,x)}\geq \gn{\bar\gb}\geq \gb$. Hence, the quantifier on $\gb$ in equation \refeq{itref} can be bounded by $x$. Thus, some elementary formula $\rho(\ga,x)$ satisfying \refeq{itref} can be constructed by the fixed point lemma.

The parametrized family of theories numerated by $\rho(\ga,x)$ will be denoted $R^\ga(\gs)$ and the formula $\rho(\ga,x)$ will be more suggestively written as $x\in R^\ga(\gs)$. Then, equation \refeq{itref} can be interpreted as saying that
$R^0(\gs) =_\EA \gs $ and, if $\ga\succ 0$,
$$
R^\ga(\gs) =_\EA \textstyle{\bigcup}\{R(R^\gb(\gs)):\gb\prec \ga\}.
$$

\bl\ \label{monit} Suppose $R$ is uniformly definable.
\benr
\item
If $0\prec \ga\preceq\gb$ then $R^\gb(\gs)\leq_\EA R^\ga(\gs)$;
\item \emph{$\EA\vdash \al{\ga,\gb}(0\prec \ga\prec\gb \to \bb R^\gb(\gs)\leq R^\ga(\gs)\qq).$}
\eenr
\el
\bp\ Obviously, Claim (i) follows from Claim (ii). For the latter we unwind the definition of $\rho(\ga,x)$ and prove within $\EA$
\beq \label{monrho} \al{\ga,\gb}(0\prec \ga\prec\gb \to \al{x}(\rho(\ga,x)\to \rho(\gb,x)).\eeq
This is sufficient to obtain from the same premise $\al{x}(\Box_{\rho(\ga,\cdot)}(x)\to \Box_{\rho(\gb,\cdot)}(x)).$

For a proof of \refeq{monrho} we reason within $\EA$: If $\rho(\ga,x)$ and $\ga\not\approx 0$ then there is a $\gy\prec\ga$ such that $\Ax_R(x,\gn{\rho(\ol \gy, x)})$. By the provable transitivity of $\prec$ from $\ga\prec\gb$ we obtain $\gy\prec\gb$, hence $\rho(\gb,x)$, q.e.d. \ep

\bl\ \label{mon2} Suppose $R$ is reflexively monotone.
If $\tau\leq_\EA\gs$ then $R^\ga(\tau)\leq_\EA R^\ga(\gs)$ and, moreover, \emph{$\EA\vdash \al{\ga} \bb R^\ga(\tau)\leq R^\ga(\gs)\qq.$}
\el

\bp\ We argue by reflexive induction similarly to \cite{Bek99b}, that is, we prove in $\EA$ that
$$\al{\gb\prec\ga}\Box_\EA\al{x}(\Box_{R^{\ol{\gb}}(\gs)}(x)\to \Box_{R^{\ol{\gb}}(\tau)}(x)) \to \al{x}(\Box_{R^\ga(\gs)}(x)\to \Box_{R^\ga(\tau)}(x))$$ and then apply L\"ob's theorem in $\EA$. Assume $\tau\leq_\EA\gs$.

Reason within $\EA$: If $\Box_{R^\ga(\gs)}(x)$ then either $\ga\approx 0\land \Box_\gs(x)$, or there is a $\gb\prec\ga$ such that $\Box_{R(R^{\ol{\gb}}(\gs))}(x)$. In the first case we obtain $\Box_{\tau}(x)$ by the external assumption $\tau\leq_\EA\gs$ and are done. In the second case, by the premise and the reflexive monotonicity of $R$ we obtain $\Box_{R(R^{\ol{\gb}}(\tau))}(x)$ which yields $\Box_{R^\ga(\tau)}(x)$.
\ep

\bcor The iteration of $R$ along $(D,\prec)$ is uniquely defined, that is, equation \emph{\refeq{itref}} has a unique solution modulo $=_\EA$.
\ecor

\bl \label{semiid}
Suppose $R$ is reflexively monotone and semi-idempotent. Then
\benr
\item
If $0\prec \ga$ then $R(R^\ga(\gs))\leq_\EA R^\ga(\gs)$;
\item \emph{$\EA\vdash \al{\ga}(0\prec \ga \to \bb R(R^\ga(\gs))\leq R^\ga(\gs)\qq).$}
\eenr
\el

\bp\ Claim (i) follows from (ii). For the latter, it is sufficient to prove the claim within $\EA+\BS_1$ and refer to the $\Pi_2^0$-conservativity of $\BS_1$ over $\EA$~(cf \cite{HP}).

Reason in $\EA+\BS_1$: If $0\prec\ga$ and $x\in R^\ga(\gs)$, then there is a $\gb\prec\ga$ such that $x\in R(R^\gb(\gs))$. We consider two cases. If $0\prec \gb$ then, since (provably) $\gb\prec\ga$, by Lemma~\ref{monit} we have $R^\ga(\gs)\leq_\EA R^\gb(\gs)$. By the reflexive monotonicity of $R$ we obtain $R(R^\ga(\gs))\leq R(R^\gb(\gs))$. Hence, $\Box_{R(R^\ga(\gs))}(x)$ and we are done.

If $\gb\approx 0$ then by the definition $R^\gb(\gs)=_\EA \gs$. Hence, by the reflexive monotonicity of $R$, $R(\gs)\leq R(R^\gb(\gs))$.
Since $0\prec\ga$, by the definition $R^\ga(\gs)\leq_\EA R(R^0(\gs))\leq_\EA R(\gs)$. It follows that $R(R^\ga(\gs))\leq R(R(\gs))\leq R(\gs)$ and therefore $\Box_{R(R^\ga(\gs))}(x)$.

Thus, using $\BS_1$ we may conclude that $0\prec\ga$ implies
$$\al{x}(\Box_{R^\ga(\gs)}(x)\to \Box_{R(R^\ga(\gs))}(x)),$$
as required. \ep

The following lemma is most naturally stated for elementary pre-wellorderings equipped with elementary formulas $\Suc(\ga,\gb)$ expressing ``$\gb$ is a successor of $\ga$'' and
$\Lim(\ga)$ expressing ``$\ga$ is a limit'' that provably in $\EA$ satisfy their defining properties:
\begin{eqnarray*}
\al{\ga,\gb}(\Suc(\ga,\gb)& \eqv & (\ga\prec\gb \land \al{\gy}(\gy\prec\gb\to \gy\preceq \ga)); \\
\al{\ga} (\Lim(\ga) & \eqv & \neg\:\ga\approx 0\land \al{\gb}(\gb\prec\ga\to \ex{\gy}(\gb\prec\gy \land \gy\prec\ga))).
\end{eqnarray*}

\bl Suppose $R$ is reflexively monotone and semi-idempotent. Then
\benr
\item $R^\ga(\gs)=_\EA\gs$ if $\ga\approx 0$, \item $R^{\gb}(\gs)=_\EA R(R^\ga(\gs))$ if $\Suc(\ga,\gb)$,  \item $R^{\gl}(\gs)=_\EA\bigwedge\{R^\ga(\gs):0\prec\ga\prec\gl\}$ if $\Lim(\gl)$.
\eenr
\el
Here $\bigwedge\{R^\ga(\gs):0\prec\ga\prec\gl\}$ denotes the G\"odelian theory numerated by $$ \ex{\ga}(0\prec\ga\prec\underline{\gl} \land x\in R^\ga(\gs)).$$
\bp\ Claim (i) is easy. For  Claim (ii) assume $\Suc(\ga,\gb)$. The implication $R^{\gb}(\gs)\leq_\EA R(R^\ga(\gs))$ is easy, since $\ga\prec\gb$ and this fact is provable in $\EA$. For the opposite implication it is sufficient to prove in $\EA+\BS_1$: $$\al{x}(x\in R^{\gb}(\gs)\to \Box_{R(R^\ga(\gs))}(x)).$$ Then one will be able to conclude using $\BS_1$ that
$\al{x}(\Box_{R^{\gb}(\gs)}(x) \to \Box_{R(R^\ga(\gs))}(x))$ and then appeal to the $\Pi_2^0$-conservativity of $\BS_1$ over $\EA$.

Reason in $\EA+\BS_1$: Assume $x\in R^{\gb}(\gs)$ then (since $\gb\not\approx 0$) there is a $\gy\prec \gb$ such that  $x\in R(R^\gy(\gs))$. If $\gy\approx\ga$ then $x\in R(R^\ga(\gs))$ and we are done. Otherwise, $\gy\prec \ga$ and one has $R^\ga(\gs)\leq R(R^\gy(\gs))$. On the other hand, by Lemma~\ref{semiid}, $R(R^\ga(\gs))\leq R^\ga(\gs)$. Hence, $R(R^\ga(\gs))\leq R(R^\gy(\gs))$, therefore $\Box_{R(R^\ga(\gs))}(x)$ as required.

\ignore{
If $\gy\succ 0$ then we have $R^\ga(\gs)\leq_\EA R^\gy(\gs)$ by Lemma \ref{monit} (ii). By the reflexive monotonicity of $R$ we obtain $R(R^\ga(\gs))\leq R(R^\gy(\gs))$, whence $\Box_{R(R^\ga(\gs))}(x)$.

If $\gy\approx 0$ then $R^\gy(\gs)\leq_\EA \gs$. Hence, by reflexive monotonicity $R(R^\gy(\gs))\leq R(\gs)$. Since $\gy\prec\ga$ we also have $R^\ga(\gs)\leq R(R^\gy(\gs))$.
$x\in R(\gs)$. We have $R^\ga(\gs)\leq_\EA R(\gs)$ since $\ga\succ 0$. Hence, $R(R^\ga(\gs))\leq_\EA R(R(\gs))\leq_\EA R(\gs)$ by the semi-idempotence of $R$. So, from $x\in R(\gs)$ we infer $\Box_{R(R^\ga(\gs))}(x)$.
}

To prove Claim (iii) we argue in a similar manner. Assume $\Lim(\gl)$, then this fact is also provable in $\EA$. To prove the implication from left to right we reason in $\EA+\BS_1$:

Assume $x\in R^\gl(\gs)$. Since $\gl\not\approx 0$ there is a $\gb\prec\gl$ such that $x\in R(R^\gb(\gs))$. Since $\Lim(\gl)$ there is an $\ga$ such that $\gb\prec\ga\prec\gl$. Then $R^\ga(\gs)\leq_\EA R(R^\gb(\gs))$ and hence $\Box_{R^\ga(\gs)}(x)$.

From right to left we reason in $\EA+\BS_1$:
Assume $0\prec\ga\prec\gl$ and $x\in R^\ga(\gs)$. Since $\ga\prec\gl$ we have by definition $R^\gl(\gs)\leq R(R^\ga(\gs))$.
On the other hand,  since $\ga\succ 0$, by Lemma \ref{semiid} we have $R(R^\ga(\gs))\leq R^\ga(\gs)$. Then $R^\gl(\gs)\leq R^\ga(\gs)$ and $\Box_{R^\gl(\gs)}(x)$, as required.
\ep


\section{Expressibility of iterated reflection}

In this section we confuse the arithmetical and reflection calculus notation. We write $\Diamond_n$ for $\mR_n$ and $\nab_n$ for $\iPi_{n+1}$.
Our goal is to show that iterated operators $\Diamond_n^\ga$, for natural ordinal notations $\ga<\ge_0$, are expressible in the language of $\Rcn$. We will rely on the so-called \emph{reduction property} (cf.~\cite{Bek04}, the present version is somewhat more general and follows from \cite[Theorem 2]{Bek99b}, see also \cite{Bek17}).

 Let $\EA^+$ denote the theory $\mR_1(\EA)$ which is known to be equivalent to $\EA+\text{Supexp}$. Theories in this and the following section will be G\"odelian extensions of $\EA^+$. We could have worked more generally over $\EA$ at the cost of replacing the reflection and conservativity operators of $\fG_\EA$ by their analogs stated for cut-free provability (see~\cite[Appendix C]{Bek99b}). Taking cut-free version of $\EA$ as our base G\"odelian theory seems to be a better choice for proof-theoretic applications. However, for simplicity we prefer to strengthen our base theory to $\EA^+$ as it was done in some previous papers we would like to refer to.

 Working in $\fG_{\EA^+}$ we write $\Diamond_{n,\gs}(\tau)$ for $\Diamond_n(\gs\land\tau)$. Obviously,  $\Diamond_{n,\gs}$ is a monotone semi-idempotent operator, for each $\gs$. Also, $1$ will stand for $1_{\EA^+}$.

\bt[reduction property] \label{redprop}
For all $\gs\in\fG_{\EA^+}$, $n\in \gw$, $$\Diamond_{n,\gs}^\gw(1)=_{\EA^+} \nab_n\Diamond_{n+1}(\gs).$$
\et

We also remark that the theory $\Diamond_{n,\gs}^\gw(1)$ is equivalent to the one  axiomatized over $\EA$ by the union of theories $\{\mQ_n^k(\gs):k<\gw\}$, where $\mQ_n^0(\gs):=\Diamond_n\gs$ and $\mQ_n^{k+1}(\gs):=\Diamond_n(\gs\land \mQ_n^k(\gs))$ are defined by formulas in one variable of $\Rc$. The corresponding G\"odelian theory taken with its natural numeration will be denoted $\bigwedge_{k<\gw}\mQ_n^k(\gs)$.

Concerning these formulas we note three well-known facts.

\bl\ \label{Q}
Provably in $\EA$,
\ben
\item $\al{B\in \Wo_n} \al{k} \mQ_n^{k+1}(B)\vdash_\Rc \mQ_n^{k}(B)\land \Diamond_n \mQ_n^{k}(B)$;
\item $\al{B\in \Wo_n} \al{k} \mQ_n^{k}(B)<_n \Diamond_{n+1} B$;
\item $\al{B\in \Wo_n} \al{k} \ex{A\in \Wo_n} \mQ_n^{k}(B)=_\Rc A$.
\een
\el

The first two of these claims are proved by an easy induction on $k$. The third one is a consequence of a more general theorem that any variable-free formula of $\Rc$ is equivalent to a word. An explicit rule for calculating such an $A$ is also well-known and related to the so-called \emph{Worm sequence}, see \cite[Lemma 5.9]{Bek05en}.

We consider the set of words $(\Wo_n,<_n)$ modulo equivalence in $\Rc$, together with its natural representation in $\EA$, as an elementary pre-wellordering. Recall that, for each $A\in \Wo_n$,  $o_n(A)$ denotes the order type of $\{B<_n A:B\in \Wo_n\}$ modulo $=_\Rc$. In a formalized context, the ordinal $o_n(A)$ is represented by its notation, the word $A$, however we still write $o_n(A)$, as it reminds us that $A$ must be viewed as an ordinal and indicates which system of ordinal notation is considered.
From the reduction property we obtain the following theorem that was stated as Theorem 6 in \cite{Bek04} in a somewhat different way. We provide a proof for the reader's convenience, though it is nearly the same as in \cite{Bek04}.

For a word $A$ and a G\"odelian theory $\gs\in \fG_{\EA^+}$, let $A^*(\gs)$ denote the interpretation of the formula $A[p/\top]$ in $\fG_{\EA^+}$ sending $p$ to $\gs$.

\bt \label{itref1} For all words $A\in \Wo_n\setminus\{\top\}$, in $\fG_{\EA^+}$ there holds $$\nab_n A^*(\gs) =_{\EA^+} \Diamond_{n}^{o_n(A)}(\gs).$$
\et

\bp\ We argue by reflexive induction in $\EA^+$ and prove that, for all $\gs\in\fG_{\EA^+}$ and all $n<\gw$,
$$\EA^+\vdash \al{B<_n A}\bb \nab_n B^*(\gs) =_{\EA^+} \Diamond_{n}^{o_n(B)}(\gs)\qq \to \bb\nab_n A^*(\gs) = \Diamond_{n}^{o_n(A)}(\gs)\qq.$$

Arguing inside $\EA^+$, we will omit the quotation marks and read the expressions $\tau\leq\nu$ as  $\al{x}(\Box_\nu(x)\to \Box_\tau(x))$ and $\tau=\nu$ as $\al{x}(\Box_\tau(x)\eqv \Box_\nu(x))$.


If $A\circeq\Diamond_n B$ then $o_n(A)=o_n(B)+1$. If $B\circeq\top$ the claim follows since $\nab_n A^*(\gs)=\nab_n\Diamond_n\gs=\Diamond_n(\gs).$
If $B\not\circeq \top$ we have by the reflexive induction hypothesis $\nab_n B^*(\gs) =_{\EA^+} \Diamond_{n}^{o_n(B)}(\gs)$. It follows that $\Diamond_{n}(\Diamond_{n}^{o_n(B)}(\gs))= \Diamond_{n}\nab_n B^*(\gs)= \Diamond_{n} B^*(\gs).$ Therefore, we obtain $$\Diamond_{n}^{o_n(A)}(\gs)=\Diamond_{n}
(\Diamond_{n}^{o_n(B)}(\gs))= \Diamond_{n} B^*(\gs)
= A^*(\gs)= \nab_n A^*(\gs).$$

If $A\circeq\Diamond_{m+1} B$ with $m\geq n$ then $\nab_n A^*(\gs)=\nab_n\nab_{m}\Diamond_{m+1} B^*(\gs)$. By the reduction property $\nab_{m}\Diamond_{m+1} B^*(\gs)= \bigwedge_{k<\gw} \mQ_m^k(B^*(\gs))$. Moreover, by Lemma~\ref{Q}~(i), if a sentence is provable in $\bigwedge_{k<\gw} \mQ_m^k(B^*(\gs))$, it must be provable in $ \mQ_m^k(B^*(\gs))$, for some $k<\gw$.
Hence, we can infer
$$\textstyle\nab_n A^*(\gs)= \nab_n\bigwedge_{k<\gw} (\mQ_m^k(B^*(\gs))=\bigwedge_{k<\gw} \nab_n(\mQ_m^k(B^*(\gs)) =\bigwedge_{k<\gw} \Diamond_{n}\mQ_m^k(B^*(\gs)).$$
By Lemma \ref{Q}(ii) and (iii), each of $\mQ_m^k(B)$ is $<_n$-below $A\circeq\Diamond_{m+1} B$ and is equivalent to a word in $\Wo_n$. Hence,
$$\textstyle \bigwedge_{C<_n A} \Diamond_{n} C^*(\gs)\leq \bigwedge_{k<\gw} \Diamond_{n}(\mQ_m^k(B^*(\gs)).$$
By the reflexive induction hypothesis, for each $C<_n A$ we have $$\Diamond_{n}C^*(\gs) =\Diamond_{n}\nab_n C^*(\gs) = \Diamond_{n} \Diamond_{n}^{o_n(C)}(\gs).$$ (If $C=\top$ the claim holds trivially.)
It follows that
\begin{multline*}\textstyle\Diamond_{n}^{o_n(A)}(\gs)= \bigwedge_{C<_n A} \Diamond_{n}\Diamond_{n}^{o_n(C)}(\gs) = \bigwedge_{C<_n A}\Diamond_{n} C^*(\gs) \leq \\ \textstyle\leq \bigwedge_{k<\gw} \Diamond_{n}\mQ_m^k(B^*(\gs)) = \nab_n A^*(\gs).
\end{multline*}

On the other hand, if $C<_n A$ then $A^*(\gs)\leq \Diamond_{n} C^*(\gs)$ and $\nab_n A^*(\gs)\leq \nab_n\Diamond_{n} C^*(\gs)\leq \Diamond_{n} C^*(\gs)$. Hence, $$\textstyle\nab_n A^*(\gs)\leq  \bigwedge_{C<_n A}\Diamond_{n} C^*(\gs) =  \Diamond_{n}^{o_n(A)}(\gs).$$
Thus, we have proved $\nab_n A^*(\gs) = \Diamond_{n}^{o_n(A)}(\gs)$, as required. \ep

For ordinals $\ga<\ge_0$, let $\mA^n_\ga\in \Wo_n$ denote a canonical notation for $\ga$ in the system $(\Wo_n,<_n)$. Thus, $o_n(\mA^n_\ga)= \ga$.
We are going to show that the operations $\Diamond_n^\ga$ are expressible in $\Rcn$ in the following sense.

\bt \label{expr-it} For each $n<\gw$ and $0<\ga<\ge_0$ there is an $\Rc$-formula $A(p)$ such that $\al{\gs\in\fG_{\EA^+}} \Diamond_n^\ga(\gs) =_{\EA^+} \nab_n A^*(\gs)$.
\et

\bp\ Take $ A(p):=\mA^n_\ga[p/\top]$ and apply Theorem \ref{itref1}. \ep

\section{Proof-theoretic $\Pi_{n+1}^0$-ordinals and conservativity spectra} \label{cns-sp}

\renewcommand{\ord}{\mathrm{ord}}

Let $S$ be a G\"odelian extension of $\EA^+$ and $(\Omega,<)$ a fixed  elementary recursive well-ordering. In this section we additionally assume that $\Omega$ is an epsilon number and is equipped with elementary terms representing the ordinal constants and functions $0,1,+,\cdot,\gw^x$. These functions should provably in $\EA$ satisfy some minimal natural axioms NWO listed in \cite{Bek95}. We call such well-orderings \emph{nice}. Recall the following definitions from \cite{Bek04} (writing $1$ for $1_{\EA^+}$):

\bi
\item \emph{$\Pi_{n+1}^0$-ordinal of $S$}, denoted $\ord_n(S)$, is the supremum of all $\ga\in\Omega$ such that $S\vdash \mR_n^\ga(1)$;
\item $S$ is \emph{$\Pi_{n+1}^0$-regular} if $S$ is $\Pi_{n+1}^0$-conservative over $\mR_n^\ga(1)$, for some $\ga\in\Omega$.
\ei

The following basic proposition states that $\Pi^0_{n+1}$-ordinals are insensitive to  $\Pi^0_{n+1}$-conservative extensions and to extensions by consistent $\Sigma^0_{n+1}$-axioms.

\bpr \label{insense} For any $S,T$ and a nice well-ordering $\Omega$, for all $n\in\gw$,
\benr\item If $S\vdash \Pi_{n+1}(T)$ then $\ord_n(S)\geq \ord_n(T)$;
\item If $T$ is axiomatized by $\Sigma_{n+1}^0$-sentences and $S\cup T$ is consistent, then \\ $\ord_n(S\cup T)=\ord_n(S)$.
    \eenr
\epr
\bp\ The first claim follows from the fact that $\mR_n^\ga(1)$ is a $\Pi_{n+1}^0$-axiomatized theory. The second claim follows from the well-known result by Kreisel and L\'evy~\cite{KrL} that $\mR_n(U)$ is not contained in any consistent $\Sigma_{n+1}^0$-axiomatized extension of $U$. \ep\

We refer the readers to \cite{Bek04} or \cite{Bek15a} for an extended discussion of proof-theoretic $\Pi_{n+1}^0$-ordinals. In this paper we consider the sequences of $\Pi_{n+1}^0$-ordinals associated with a given system. Such sequences as objects of study first appeared in the work of Joost Joosten~\cite{Joo15a}. He showed for theories between $\EA^+$ and $\PA$ that their conservativity spectra correspond to decreasing sequences of ordinals below $\ge_0$ of a certain kind, that is, to the points in the so-called Ignatiev frame. We reproduce this interesting characterization here in a slightly more general and streamlined way and also show its tight relationship with the fat normal forms for $\Rcn$.

\bd \emph{Conservativity spectrum of $S$} is the sequence $(\ga_0,\ga_1,\ga_2,\dots)$ such that $\ga_i=\ord_i(S)$.
\ed

Here are some examples of theories and their spectra (the results are either well-known and/or can be found in \cite{Bek99b}):
\ben \item
$I\Sigma_1: \ (\gw^\gw,\gw,1,0,0,\dots)$;\qquad $\PRA: \ (\gw^\gw,\gw,0,0,0,\dots)$;
\item $\PA: \ (\ge_0,\ge_0,\ge_0,\dots )$; \qquad $\PA+\Con(\PA): \ (\ge_0\cdot 2,\ge_0,\ge_0,\dots )$ \item
$\PA+\mR_1(\PA):  \quad (\ge_0^2,\ge_0\cdot 2,\ge_0,\ge_0,\dots )$.
\een

We will need the following auxiliary lemma concerning the iterations of the reflection operators $\mR_n$ on $\fG_{T}$, for any G\"odelian extension $T$ of $\EA$.
\bl \label{it-j}
Let $(D,\preceq,0)$ be an elementary pre-wellordering, then for all $\gs_1,\gs_2\in \fG_T$ there holds
$$\al{\ga\succ 0}\ \mR_{n+1}^\ga(\gs_1)\land_T \mR_n(\gs_2)=_T \mR_{n+1}^\ga(\gs_1\land_T \mR_n(\gs_2)).$$
\el

\bp\ The proof is routine by reflexive induction on $\ga$ using the $\Rc$-identity
$$\mR_{n+1}(\gs_1)\land_T \mR_n(\gs_2)=_T \mR_{n+1}(\gs_1\land_T \mR_n(\gs_2)).$$
\ep

The following proposition provides a necessary condition for a sequence of ordinals to be a conservativity spectrum.

\bpr\ \label{spec-ign} For any $S$ and a nice well-ordering $\Omega$, for all $n\in\gw$,
\benr \item $\ord_{n+1}(S)\leq \ell(\ord_n(S))$;
\item If $S$ is $\Pi_{n+1}^0$-regular and $n>0$, then $\al{i<n}\ord_{i}(S) = \gw^{\ord_{i+1}(S)}$. \eenr\epr
\bp\
For (i) let $\vec\ga$ denote the conservativity spectrum of $S$ and assume $\ga_{n+1}>\ell(\ga_n)$. Select a $\gy$ such that $\ga_n=\gy+\gw^{\ell(\ga_n)}$. Notice that $$S\vdash \mR_{n+1}^{\ga_{n+1}}(1)\land_{\EA^+} \mR_n^{\gy+1}(1).$$ Then by
Lemma \ref{it-j} we obtain
$$\mR_{n+1}^{\ga_{n+1}}(1)\land_{\EA^+} \mR_n^{\gy+1}(1)=_{\EA^+}  \mR_{n+1}^{\ga_{n+1}}(\mR_n^{\gy+1}(1)).$$
By Theorem 3 of \cite{Bek04} we have: $$\al{\gb\in\Omega}\ \iPi_{n+1} \mR_{n+1}^\gb(\gs) =_{\EA^+} \mR_n^{\gw^\gb}(\gs),$$ for all $\Pi_{n+1}$-axiomatized extensions $\gs$ of $\EA^+$. Hence, $$S\vdash \mR_{n+1}^{\ga_{n+1}}(\mR_n^{\gy+1}(1))\vdash \mR_n^{\gy+\gw^{\ga_{n+1}}}(1).$$
It follows that $\ga_n=\ord_n(S)\geq \gy+\gw^{\ga_{n+1}}$. On the other hand, by our assumption $\gy+\gw^{\ga_{n+1}}>\gy+\gw^{\ell(\ga_n)}=\ga_n$, a contradiction.

Since a $\Pi_m$-regular theory is $\Pi_i$-regular, for all $i<m$, it is sufficient to prove the claim for $i=n-1$. If $S$ is $\Pi_{n+1}^0$-regular, then $\iPi_{n+1}(S)= \mR_{n}^\ga(1)$ where $\ga=\ord_n(S)$. It follows that $\iPi_{n}(S)=\iPi_n(\iPi_{n+1}(S))= \iPi_n(\mR_{n}^\ga(1))=\mR_{n-1}^{\gw^\ga}(1)$.
\ep

We consider the ordering of words $(\Wo_i,<_i)$ modulo equivalence in $\Rc$ as a nice well-ordering of length $\ge_0$. As before, the order type of a word $A_i$ within $(\Wo_i,<_i)$ is denoted $o_i(A_i)$. A direct correspondence between fat normal forms and conservativity spectra is expressed by the following theorem.

\bt \label{spect} Let $A\circeq\nab_0A_0\land\nab_{1}A_{1}\land\dots\land\nab_{k}A_{k},$ for some $k$, where $A_n\in \Wo_n$ for all $n\leq k$, be in the fat normal form. Let $A^*$ denote the interpretation of $A$ in $\fG_{\EA^+}$.
Then, $A_n$ represents the $\Pi_{n+1}^0$-ordinal of $A^*$: $o_n(A_n)=\ord_n(A^*).$ Moreover, $A^*$ is equivalent to the union of progressions, that is, in $\fG_{\EA^+}$
\beq \label{utf} A^*=_{\EA^+} \mR_0^{o_0(A_0)}(1)\land \mR_1^{o_1(A_1)}(1)\land\dots\land \mR_k^{o_k(A_k)}(1).\eeq
\et

\bp\ Firstly, by applying Proposition \ref{insense} we observe that $$\ord_n(A^*)=\ord_n((\nab_nA_n\land\nab_{n+1}A_{n+1}\land\dots\land\nab_{k}A_{k})^*)= \ord_n(A_n^*).$$ The first equality holds, because the deleted part of the normal form of $A$ is interpreted as a true $\Pi^0_{n}$-theory. The second equality holds, since the remaining part of the fat normal form of $A$ is $\Pi_{n+1}^0$-conservative over $A_n^*$:
$$\nab_nA_n=_\Rcn \nab_n(\nab_nA_n\land\nab_{n+1}A_{n+1}\land\dots\land\nab_{k}A_{k}).$$
Then, by Theorem~\ref{itref1},
\begin{equation}
\Pi_{n+1} (A_n^*) =_{\EA^+} \mR_n^{o_n(A_n)}(1), \label{tte}
\end{equation} hence $A_n^*$ is $\Pi_{n+1}^0$-regular and $\ord_n(A_n^*)=o_n(A_n)$. Moreover, equation \refeq{tte} also  yields representation \refeq{utf} of $A^*$ as a union of progressions. \ep

Joosten~\cite{Joo15a} calls the representations of theories as the unions of Turing progressions \emph{Turing--Taylor expansions}. Thus, the fat normal form of $A$ represents the Turing--Taylor expansion of $A^*$ by way of \refeq{utf}.

Notice that Theorem \ref{spect} also yields another way of showing that the fat normal form of $A$ is unique. We will come back to the topic of conservativity spectra after we discuss the Ignatiev frame and its associated $\Rcn$-algebra.

\section{Ignatiev frame and Ignatiev $\Rcn$-algebra}

In this and the following section we characterize in several ways the Lindenbaum--Tarski algebra of the variable-free fragment of $\Rcn$. It turns out that this structure is tightly related to the so-called \emph{Ignatiev's} Kripke frame. This frame, denoted here $\Ig$, has been introduced by Konstantin Ignatiev~\cite{Ign93} as a universal frame for the variable-free fragment of Japaridze's logic GLP. Later this frame has been slightly modified and studied in more detail in \cite{BJV,Ica09}. In particular, Thomas Icard established a detailed relationship between $\Ig$ and the canonical frame for the variable-free fragment of $\Glp$ and used it to define a complete  topological semantics for this fragment. David Fern\'andez and Joost Joosten \cite{FJ13a} generalized $\Ig$ to a version of GLP with transfinitely many modalities. Ignatiev's frame is defined constructively (`coordinatewise') as follows.

Let $\oI$ denote the set of all $\gw$-sequences of ordinals $\vec \ga=(\ga_0,\ga_1,\dots)$ such that $\ga_i\leq \ge_0$ and $\ga_{i+1}\leq \ell(\ga_i)$, for all $i\in\gw$. Here, the function $\ell$ is defined by: $\ell(\gb)=0$ if $\gb=0$, and $\ell(\gb)=\gy$ if $\gb=\gd+\gw^{\gy}$, for some $\gd,\gy$. Thus, all sequences of $\oI$, with the exception of identically $\ge_0$, are eventually zero. Elements of $\oI$ will also be called \emph{$\ell$-sequences}.

Relations $R_n$ on $\oI$ are defined by:
$$\vec\ga R_n\vec\gb \iff (\al{i<n} \ga_i=\gb_i \text{ and } \ga_n>\gb_n).$$
The structure $\oIg=(\oI,(R_n)_{n\in\gw})$ is called \emph{the extended Ignatiev frame} (see \cite{Ica09}).
The Ignatiev frame is its restriction to the subset $\mI$ of all sequences $\vec\ga\in\oI$ such that $\al{i\in\gw} \ga_i<\ge_0$. This subset is upwards closed w.r.t.\ all relations $R_n$, hence the evaluation of the variable-free $\Rc$-formulas (and $\Glp$-formulas) in $\Ig$ and in $\oIg$ coincide. We denote by $\Ig,\vec\ga\fc \phi$ the truth of a $\Glp$-formula $\phi$ at a node $\vec\ga$ of $\Ig$. If $\phi$ is variable-free, the set $\{\vec\ga\in \mI: \Ig,\vec\ga\fc \phi\}$ will be denoted $v(\phi)$.

The following important theorem is a corollary of the results of Ignatiev but, in fact, has an easier direct proof (which we omit for the reasons of brevity).

\bpr For any variable-free formulas $A,B$ of $\Rc$, $A\vdash_\Rc B$ iff $\cI,\vec\ga\fc A\to B$, for all $\vec\ga\in \mI$.
\epr

The set of sequences $\vec\ga\in \mI$ such that $\al{i<\gw} \ga_{i+1} = \ell(\ga_i)$ is called \emph{the main axis} of $\Ig$ and is denoted $\mO$. Obviously, a sequence in $\mO$ is uniquely determined by its initial element $\ga_0$, hence $\mO$ naturally corresponds to the ordinals up to $\ge_0$.
We can also associate with every word $A\in\Wo$ an element $\iota(A)\in\mO$ by letting $$\iota(A):=(o(A),\ell(o(A)),\dots,\ell^{(n)}(o(A)),\dots).$$
The following lemma, explicitly stated by Thomas Icard \cite[Lemma 3.8]{Ica09} (see also another argument in~\cite[Lemma 10.2]{BJ18}), describes all the subsets of $\oIg$ definable by words (and hence by all variable-free s.p.\ formulas of $\Rc$).

\bl \label{Ica}
 Suppose $A\in\Wo$ and $\vec\ga=\iota(A)$. Then, for all $\vec\gb\in\oIg$, \ $\oIg,\vec\gb\fc A$ iff $\al{i\in\gw} \ga_i\leq \gb_i$.
\el

Our goal is to transform $\Ig$ into an $\Rcn$-algebra $\fI$ with the same domain $\mI$, that is, into an SLO satisfying $\Rcn$.
We consider the set $\oI$ equipped with the ordering $$\vec\ga\leq_\fI\vec\gb \iffdef \al{n\in\gw} \ga_n\geq \gb_n.$$ The structure $(\oI,\leq_\fI)$ can be seen as a subordering of the product ordering on the set of all $\gw$-sequences of ordinals $\leq\ge_0$, which we denote $\cE$.

A \emph{cone in $\cE$} is the set of points $E_{\vec{\ga}}:=\{\vec \gb\in\cE:\vec\gb\leq_\fI \vec\ga\}$, for some $\vec\ga\in\cE$. 
A sequence $\vec\ga\in\cE$ is called \emph{bounded} if $\al{i\in\gw} \ga_i<\ge_0$ and $\ga_i\neq 0$ for only finitely many $i\in\gw$. Obviously, each $\vec\ga\in \mI$ is bounded.

\bl\ \label{glb}
Suppose $\vec\ga\in\cE$ is bounded. Then $E_{\vec{\ga}}\cap \mI$ is not empty and has a  greatest point $\vec\gb$ w.r.t.\ $\leq_\fI$.

\el
\bp\ Let $n\in\gw$ be the largest number such that $\ga_n\neq 0$. Consider the sequence $\vec\gb$ such that $\gb_i=0$ for all $i>n$, $\gb_n:=\ga_n$, and, for all $i<n$:
$$\gb_i:= \begin{cases}
\ga_i, & \text{if $\ell(\ga_i)\geq \gb_{i+1}$,} \\
\ga_i+\gw^{\gb_{i+1}}, & \text{otherwise}.
\end{cases}
$$
It is easy to see that $\vec\gb$ is the greatest point of $E_{\vec\ga}\cap\mI$. Also notice that $\vec\gb$ can be effectively computed from $\vec\ga$.
\ep

\bcor $(\mI,\leq_\fI)$ is a meet-semilattice with top.
\ecor
\bp\ Let $\vec\ga,\vec\gb\in \mI$. The sequence $\vec\gy:=(\max(\ga_i,\gb_i))_{i<\gw}$ is the g.l.b.\ of $\vec\ga$ and $\vec\gb$ in $\cE$ and is bounded. By Lemma \ref{glb}, $E_{\vec\gy}\cap \mI$ has a greatest point, which has to be the g.l.b.\ of $\vec\ga$ and $\vec\gb$ in $\mI$. \ep

We denote by $\land_\fI$ the meet operation of this semilattice.
A nonempty set $C_\ga:=E_{\vec{\ga}}\cap \mI$ is called a \emph{cone in $\Ig$}. The set of all cones in $\Ig$ ordered by inclusion is denoted $\mC(\Ig)$. The orderings $(\mC(\Ig),\subseteq)$ and $(\mI,\leq_\fI)$ are isomorphic by the map $\vec\ga\mapsto C_{\vec\ga}$. So, we have

\bcor \label{C-hom}
For all $\vec\ga,\vec\gb\in \mI$, $C_{\vec\ga\land_\fI\vec\gb} =C_{\vec\ga}\cap C_{\vec\gb}$.
\ecor

Let $\mC(\mO)$ denote the set $\{C_{\vec\ga}: \vec\ga\in\mO\}$ of all cones in $\Ig$ generated by the points of the main axis. For all $X\subseteq \oI$ define $R_n^{-1}(X):=\{y\in X: \ex{x\in X} y R_n x\}$. We claim that the operations $\cap$ and $R_n^{-1}$ map cones of $\mC(\mO)$ to cones of $\mC(\mO)$. Moreover, the following proposition holds.\footnote{We do not distinguish notationally an operation on a set and its restriction to a subset.}

\bpr \label{iso1} The algebra $\fC(\mO)=(\mC(\mO);\cap,\{R_n^{-1}:n\in\gw\})$ is isomorphic to the Lindenbaum--Tarski algebra $\fL^0_{\Rc}$.
\epr
\bp\ Let $v:\Fo \to \cP(\mI)$ denote the map associating with every variable-free formula $A$ of $\Rc$ the set $v(A)$ of all points where this formula is true. By the soundness and completeness of $\Rc$ w.r.t.\ the Ignatiev model we have $v(A)=v(B)$ iff $A=_\Rc B$. Moreover, by Lemma \ref{Ica} the range of $v$ consists of all the cones of $\mC(\mO)$. So, $v$ factors to a bijective map $\bar v:\fL^0_{\Rc}\to \mC(\mO)$. The operations $\cap$ and $R_n^{-1}$ correspond to the definition of truth in a Kripke model, hence $\mC(\mO)$ is closed under these operations and $\bar v$ is an isomorphism of the respective algebras.
\ep

We remark that the work of Pakhomov~\cite{Pakh15} shows that the elementary theory of the algebra $\fL^0_{\Rc}$ is undecidable. We now define the structure of an $\Rcn$-algebra on $\mI$.

\bd\ \label{alg}
For all $n\in\gw$ we define the functions $\nab_n^\fI,\Diamond_n^\fI:\mI\to \mI$. For each element $\vec\ga=(\ga_0,\ga_1,\dots,\ga_{n},\dots)\in\mI$ let:

$\nab_n^\fI(\vec \ga):=(\ga_0,\ga_1,\dots,\ga_n,0,\dots)$;

$\Diamond_n^\fI(\vec \ga):=(\gb_0,\gb_1,\dots,\gb_{n},0,\dots)$, where $\gb_{n+1}:=0$ and $\gb_i:=\ga_i+\gw^{\gb_{i+1}}$, for all $i\leq n$.

The algebra $\fI=(\mI,\land_\fI,\{\Diamond_n^\fI,\nab_n^\fI:n\in\gw\})$ is called the \emph{Ignatiev $\Rcn$-algebra}.
\ed

The definition of the operations $\Diamond_n^\fI$ is motivated by the following lemma and its corollary.

\bl \label{Diam-hom}
Suppose $\vec\ga\in \mI$ and $\vec\gb=\Diamond_n^\fI(\vec\ga)$. Then $\vec\gb\in\mO$ and
\benr
\item $C_{\vec\gb} = \bigcap_{i\leq n} R_i^{-1}(C_{\vec\ga})$;
\item If $\vec\ga\in\mO$ then $C_{\vec\gb}=R_n^{-1}(C_{\vec\ga}).$
\eenr
\el
\bp\ (i) It is easy to see that each of the sets $R_i^{-1}(C_{\vec\ga})$, for $i\leq n$, is a cone in $\Ig$ generated by the bounded sequence $(\ga_0,\dots,\ga_{i-1},\ga_i+1,0,\dots)$ from $\cE$. Hence, the intersection of these cones is a cone generated by $(\ga_0+1,\dots,\ga_{n-1}+1,\ga_n+1,0,\dots)$. Its greatest element in $\mI$ obviously coincides with $\Diamond_n^\fI(\vec\ga)$.

(ii) Clearly, $\vec\gb\in R_n^{-1}(C_{\vec\ga})$, since $\vec\gb':=(\gb_0,\gb_1,\dots,\gb_{n-1},\ga_n,\ga_{n+1},\dots)$ satisfies $\vec\gb R_n \vec\gb'$ and $\vec\gb'\leq_\fI \vec\ga$. In the opposite direction, show by downward induction on $i\leq n$ that if $\vec\gy\in R_n^{-1}(C_{\vec\ga})$ then $\gy_i\geq \gb_i$. For $i=n$ the claim is obvious. Assume $i<n$, then $\gy_i\geq \ga_i$. Since $\ell(\gy_i)\geq \gy_{i+1}\geq \gb_{i+1}$ and $\ell(\ga_i)=\ga_{i+1}<\gb_{i+1}$, we must also have $\gy_i\geq \ga_i+\gw^{\gb_{i+1}}=\gb_i$.
\ep

\bcor \label{iso2} $\fC(\mO)$ is isomorphic to the algebra $\fO=(\mO,\land_\fI,\{\Diamond_n^\fI:n\in\gw\})$.
\ecor
\bp\ Consider the bijection $c:\vec\ga\longmapsto C_{\vec\ga}$ from $\mO$ to $\mC(\mO)$. By Corollary~\ref{C-hom} this map preserves the meet, and by Lemma~\ref{Diam-hom} it preserves the diamond modalities. \ep

We summarize the previous results in the following theorem characterizing the Lindenbaum--Tarski algebra of the variable-free fragment of $\Rc$.

\bt \label{iso3} The algebras $\fL_{\Rc}^0$, $\fC(\mO)$, $\fO$ are naturally isomorphic by the following maps:
\benr \item $\bar v:\fL_{\Rc}^0 \to \fC(\mO)$;
\item $c:\fO \to \fC(\mO)$;
\item $\bar\iota:\fL_{\Rc}^0\to \fO$.\eenr
\et
Here, for any $A\in\Fo$, $\bar\iota([A]_\Rc):= \iota (A')$, where $A'\in \Wo$ is a word such that $A=_\Rc A'$. This definition is invariant, since, for any words $A',A''$, if $A'=_\Rc A''$ then $o(A')=o(A'')$ and hence $\iota(A')=\iota(A'')$.
For a proof that (iii) is an isomorphism it is sufficient to remark that $v(A)=c(\iota(A))$, for each $A\in\Wo$, by Lemma \ref{Ica}.

Our next goal is to show that $\fI$ is isomorphic to the Lindenbaum--Tarski algebra of $\Rcn$. First, we need an auxiliary lemma.

\bl \label{mainax} For every $\vec\ga\in \mI$ and $n\in\gw$, there is an $\vec\ga'\in\mO$ such that $\vec\ga'\leq_\fI \vec\ga$ and $\Diamond_n^\fI(\vec\ga)=\Diamond_n^\fI(\vec\ga')$.
\el
\bp\ Let $\ga_n':=\ga_n$, $\al{i\geq n} \ga'_{i+1}:= \ell(\ga'_i)$, and $\al{i<n} \ga'_i:=\ga_i+\gw^{\ga'_{i+1}}$. It is easy to check that $\vec\ga'$ is as required. \ep

Let $A^\fI$ denote the value of a variable-free $\Rcn$-formula $A$ in $\fI$. The following lemma shows that $\fI$ satisfies the variable-free fragment of $\Rcn$.

\bl \label{corr} For any $A,B\in\Fon$, $A\vdash_{\Rcn} B$ implies $A^\fI\leq_\fI B^\fI$.
\el
\bp\ We argue by induction on the length of $\Rcn$-derivation. In almost all the cases the proof is routine. We consider the nontrivial case of the axiom $\Diamond_n A\land \Diamond_m B\vdash \Diamond_n(A\land \Diamond_m B)$ for $m<n$. Let $\vec\ga=A^\fI$ and $\vec\gb=B^\fI$. Using Lemma~\ref{mainax} we obtain $\vec\ga',\vec\gb'\in\mO$ such that $\vec\ga'\leq_\fI\vec\ga$, $\vec\gb'\leq_\fI\vec\gb$ and $\Diamond^{\fI}_n\vec\ga=\Diamond^{\fI}_n\vec\ga'$,  $\Diamond^{\fI}_m\vec\gb=\Diamond^{\fI}_m\vec\gb'$. By Theorem \ref{iso3} the algebra $\fO$ satisfies $\Rc$, hence $$\Diamond^{\fI}_n \vec\ga'\land_{\fI} \Diamond^{\fI}_m \vec\gb'\leq_\fI \Diamond^{\fI}_n(\vec\ga'\land_{\fI} \Diamond^{\fI}_m \vec\gb').$$
Therefore, $\Diamond^{\fI}_n \vec\ga\land_{\fI} \Diamond^{\fI}_m \vec\gb\leq_\fI \Diamond^{\fI}_n(\vec\ga'\land_{\fI} \Diamond^{\fI}_m \vec\gb)\leq_\fI \Diamond^{\fI}_n(\vec\ga \land_{\fI} \Diamond^{\fI}_m \vec\gb).$
The second inequality holds by the monotonicity of $\land_{\fI}$ and $\Diamond^\fI_n$.
 \ep

\bl \label{val-fn} Suppose $A\circeq \nab_0A_0\land \nab_1 A_1\land \dots \land \nab_n A_n$ is in the fat normal form. Then $A^\fI = (o_0(A_0),o_1(A_1),\dots , o_n(A_n), 0, \dots )$.
\el

\bp\ Firstly, since each $A_i\in\Wo_i$ we obtain from Theorem \ref{iso3} that
$$(A_i)^\fI = \iota(A_i)= (\gw_i(o_i(A_i)),\gw_{i-1}(o_i(A_i)),\dots,\gw^{o_i(A_i)},o_i(A_i),\ell(o_i(A_i)),\dots), $$ where by definition $\gw_0(\ga)=\ga$ and $\gw_{k+1}(\ga)=\gw^{\gw_k(\ga)}$.
Hence,
$$(\nab_i A_i)^\fI = (\gw_i(o_i(A_i)),\gw_{i-1}(o_i(A_i)),\dots,\gw^{o_i(A_i)},o_i(A_i),0,\dots).$$

Denote $\ol{A_i}:=\nab_iA_i\land \nab_{i+1} A_{i+1}\land \dots \land \nab_n A_n$. By downwards induction on $i\leq n$ we show that $(\ol{A_i})^\fI$ equals
\beq \label{spec} (\gw_i(o_i(A_i)),\gw_{i-1}(o_i(A_i)),\dots,o_i(A_i), o_{i+1}(A_{i+1}),\dots,o_n(A_n),0,\dots).
\eeq
For $i=n$ the claim follows from the above.
Assume $i<n$ and that the claim holds for $i+1$. Since in a fat normal form $$\nab_iA_i\vdash_\Rcn \nab_i(\nab_iA_i\land \nab_{i+1}A_{i+1}),$$ by Lemma~\ref{corr} we obtain that the sequence $(\nab_iA_i)^\fI$ coordinatewise majorizes the sequence $(\nab_i(\nab_iA_i\land \nab_{i+1}A_{i+1}))^\fI$. The former has the ordinal $o_i(A_i)$ at $i$-th position, and the latter has at the same place the least ordinal $\ga$ such that $\ga\geq o_i(A_i), \gw^{o_{i+1}(A_{i+1})}$ and $\ell(\ga)\geq o_{i+1}(A_{i+1})$. Therefore, $o_i(A_i)=\ga$ and $\ell(o_i(A_i))\geq o_{i+1}(A_{i+1})$.

Now consider the sequence
$(\ol{A_i})^\fI=(\nab_iA_i\land \ol{A_{i+1}})^\fI$. By the induction hypothesis its tail coincides with that of \refeq{spec} starting from position $i+1$. Since $\ell(o_i(A_i))\geq o_{i+1}(A_{i+1})$, the ordinal $o_i(A_i)$ occurs in it on $i$-th position. Also, for each $k<i$ we have $\gw_k(o_i(A_i))\geq \gw_k(\gw^{o_{i+1}(A_{i+1})})$. It follows that the sequence $(\ol{A_i})^\fI$ coincides with \refeq{spec}. \ep

The following corollary will be useful later on.
\bcor \label{nab-inj} For any $A,B\in\Wo$ and $n\in\gw$, if $\fI\models\nab_n A= \nab_n B$ then $A=_\Rc B$.
\ecor
\bp\ Firstly, we infer: $\fI\models \nab_0 A=\nab_0\nab_n A=\nab_0\nab_n B=\nab_0 B$. By Lemma \ref{val-fn} we conclude $o(A)=o(B)$, therefore $A=_\Rc B$. \ep

\bt\ \label{I-comp}
For all $A,B\in\Fon$, $A\vdash_\Rcn B$ iff $A^\fI\leq_\fI B^\fI$.
\et
\bp\ We must only prove the `only if' part. Moreover, it is sufficient to prove it for fat normal forms $A\circeq \nab_0A_0\land \nab_1 A_1\land \dots \land \nab_n A_n$ and $B\circeq \nab_0B_0\land \nab_1 B_1\land \dots \land \nab_m B_m$. If $A^\fI\leq_\fI B^\fI$ then by Lemma \ref{val-fn} we have $n\geq m$ and $o_i(A_i)\geq o_i(B_i)$, for each $i\leq m$. Since $A_i,B_i\in\Wo_i$, this means that $A_i\vdash_\Rc \Diamond_iB_i$ or $A_i=_\Rc B_i$. In either case we can infer $\nab_iA_i\vdash_\Rcn \nab_i B_i$ for each $i\leq m$. It follows that $A\vdash_\Rcn B$. \ep

Theorem~\ref{I-comp} essentially means the following.

\bcor The Ignatiev $\Rcn$-algebra $\fI$ is isomorphic to the Lindenbaum--Tarski algebra of the variable-free fragment of $\Rcn$.
\ecor

\section{$\fI$ as the algebra of variable-free $\Rc$-theories}

Another, perhaps even more natural, view of the Ignatiev $\Rcn$-algebra is via an  interpretation of the points of $\Ig$ as variable-free $\Rc$-theories. It nicely agrees with the arithmetical interpretation in that we can also view such a theory as an  arithmetical theory (every variable-free $\Rc$-formula corresponds to an arithmetical sentence). In this section we will presuppose that the language is variable-free and will only consider variable-free formulas and theories.

A set of strictly positive formulas $T$ is called an \emph{$\Rc$-theory} if $B\in T$ whenever there are $A_1,\dots,A_n\in T$ such that $A_1\land\dots\land A_n\vdash_\Rc B$. A theory $T$ is called \emph{improper} if $T$ coincides with the set of all strictly positive formulas, otherwise it is called \emph{proper}.\footnote{We avoid the term `consistent', for even the improper theory corresponds to a consistent set of arithmetical sentences.} A theory is called \emph{bounded} if there is a strictly positive formula $A$ such that $T\subseteq \{B: A\vdash_\Rc B\}$. We will use the following basic fact.

The set $\oI$ bears a natural topology generated as a subbase by the set of all cones in $\oIg$ and their complements. By \cite[Theorem 3.12]{Ica09}, this topology coincides with the product topology of the space $\cE$ induced on $\oI$. Obviously, for each $\Rc$-formula $A$, the set $v(A)$ is clopen. Moreover, this topology is compact and totally disconnected on $\oI$, since $\oI$ is closed in $\cE$ and $\cE$ is compact by Tychonoff theorem. As a corollary we obtain the following \emph{strong completeness} result. For each $\Rc$-theory $T$ define $v(T):=\{\vec\ga\in \mI:\Ig,\vec\ga\fc T\}$.

\bpr \label{scmpl} Let $T$ be an $\Rc$-theory and $A$ an $\Rc$-formula.
\benr
\item
$T\nvdash_\Rc A$ iff there is an $\vec\ga\in \oIg$ such that $\oIg,\vec\ga\fc T$ and $\oIg,\vec\ga\nfc A$;
\item If $T$ is bounded then $T\nvdash_\Rc A$ iff there is an $\vec\ga\in \Ig$ such that $\Ig,\vec\ga\fc T$ and $\Ig,\vec\ga\nfc A$.
\eenr
\epr

\bp\ (i) The nontrivial implication is from left to right. Assume $T\nvdash_\Rc A$. There is an increasing sequence of finite theories $(T_n)_{n\in\gw}$ such that $T=\bigcup_{n\in\gw} T_n$. By the completeness of the variable-free fragment of $\Rc$ w.r.t.\ $\Ig$
each of the sets $v(T_n)\setminus v(A)$ is nonempty and clopen. By the compactness of $\oI$ there is a point $\vec\ga\in \bigcap_{n\in\gw} v(T_n)\setminus v(A)=v(T)\setminus v(A)$.

(ii) In case $T$ is bounded we have $v(T)\supseteq v(B)$, for some word $B$. There is a bounded sequence $\vec\gb\in\cE$ such that $v(T)=E_{\vec\gb}\cap \oI$: consider the pointwise supremum of the generating points of the cones $v(T_n)$ in $\oIg$, each of which is pointwise majorized by the greatest element $B^\fI$ of $v(B)$. By Lemma \ref{glb}, the set $v(T)$ has a greatest point, say $\vec\gy\in \mI$. Since $\vec\ga\in v(T)$ we have $\vec\ga\leq_\fI \vec\gy$, hence $\Ig,\vec\gy\nfc A$. \ep

For any $\Rc$-theories $T,S$ define $T\leq_\Rc S$ iff $T\supseteq S$.
The g.l.b.\ of $T$ and $S$ in this ordering, denoted $T\land_\Rc S$, is the theory generated by the union $T\cup S$. Thus, the set $\TRC$ of all bounded variable-free $\Rc$-theories is a semilattice (it is, in fact, a lattice with $T\cap S$ the l.u.b.\ of $T$ and $S$).  The set $\{A\in\Fo:\top\vdash_\Rc A\}$ corresponds to the top of this lattice and is denoted $\top_\Rc$.

For each $\vec\ga\in\oIg$ define an $\Rc$-theory $\thr{\vec\ga}:=\{A:\oIg,\vec\ga\fc A\}$. It is easy to see that $\thr{\vec\ga}$ is bounded if $\vec\ga\in \mI$ (consider the point $\vec \gb$ on the main axis of $\Ig$ such that $\vec\gb\leq_\fI\vec\ga$ and a word $B$ such that $\iota(B)=\vec\gb$).

\bpr \label{iso-th}
\benr
\item The map $\vec\ga\mapsto [\vec\ga]$ is an isomorphism between $(\Ig,\leq_\fI)$ and the ordered set $\TRC$ of bounded $\Rc$-theories.
\item The map $v$ is an isomorphism between $\TRC$ and the ordered set $(\mC(\Ig),\subseteq)$ of cones in $\Ig$.
\eenr
\epr

\bp\
It is sufficient to prove that
\begin{enuma}
\item
The maps $\vec\ga\mapsto [\vec\ga]$ and $T\mapsto v(T)$ are order-preserving;
\item $\al{\vec\ga\in \mI} v([\vec\ga])=C_{\vec\ga}$;
\item If $v(T)=C_{\vec\ga}$ then $T=[\vec\ga]$.
\end{enuma}
Item (a) is obvious. For (b) we observe: $$\vec\gb\in v([\vec\ga]) \iff \al{A} (\Ig,\vec\ga\fc A \Imp \Ig,\vec\gb\fc A).$$
The right hand side is equivalent to $\vec\gb\leq_\fI \vec\ga$: If $\vec\gb\leq_\fI \vec\ga$ and $\Ig,\vec\ga\fc A$ then $\Ig,\vec\gb\fc A$ by Proposition \ref{Ica}. If $\vec\gb\nleqslant_\fI \vec\ga$ then there is a word $A$ such that $\Ig,\vec\ga\fc A$ and $\Ig,\vec\gb\nfc A$, by \cite[Corollary 3.9]{Ica09}. Hence, $\vec\gb\in v([\vec\ga])$ iff $\vec\gb\in C_{\vec\ga}$.

For (c) we use Proposition \ref{scmpl}. Suppose $\vec\ga\in \mI$ and $v(T)=C_{\vec\ga}$. Then $\Ig,\vec\ga\fc T$ and thus $T\subseteq [\vec\ga]$. For the opposite inclusion assume $A\in [\vec\ga]$ and $A\notin T$. By Proposition \ref{scmpl} there is a node $\vec\gb\in \mI$ such that $\Ig,\vec\gb\fc T$ and $\Ig,\vec\gb\nfc A$. Thus, $\vec\gb\in v(T)$ and, since $v(A)$ is downwards closed, $\vec\gb\nleqslant_\fI\vec\ga$. It follows that $v(T)\nsubseteq C_{\vec\ga}$.
\ep

The operations of the Ignatiev $\Rcn$-algebra can be interpreted in terms of the semilattice of bounded theories as follows. For each $T\in\TRC$ let
$\nab_n^\Rc T$ denote the $\Rc$-theory axiomatized by  $\{ \Diamond_m A: \Diamond_m A\in T \text{ and } m\leq n\}.$

\bl\ \label{iso-nab-th} For all $\vec\ga\in \fI$, $\nab^{\Rc}_n([\vec\ga])=[\nab^\fI_n\vec\ga]$.
\el

\bp\ For the inclusion $(\subseteq)$ we need to show: if $m\leq n$ and $\Diamond_mA\in [\vec\ga]$ then $\Diamond_mA\in [\nab^\fI_n\vec\ga]$. If $\Diamond_mA\in [\vec\ga]$ then $\Ig,\vec\ga\fc \Diamond_m A$, hence there is a $\vec\gb$ such that $\vec\ga R_m\vec\gb$ and $\Ig,\vec\gb\fc A$. So, we have $\al{i<m}\, \ga_i=\gb_i$ and $\ga_m>\gb_m$. Since $m\leq n$, the node $\nab^\fI_n\vec\ga$ has the same coordinates as $\vec\ga$ for all $i\leq m$. Therefore, $(\nab^\fI_n\vec\ga) R_m \vec\gb$ and $\Ig,(\nab^\fI_n\vec\ga)\fc \Diamond_m A$.

For the inclusion $(\supseteq)$ we consider any node $\vec\gy\in \mI$ such that $\Ig,\vec\gy\fc \nab^{\Rc}_n[\vec\ga]$ and show that $\Ig,\vec\gy\fc [\nab^\fI_n\vec\ga]$. This means that $v(\nab^{\Rc}_n[\vec\ga])\subseteq v([\nab^\fI_n\vec\ga])$ and hence $\nab^{\Rc}_n([\vec\ga])\supseteq [\nab^\fI_n\vec\ga]$ by Proposition~\ref{iso-th}.

Assume $\Ig,\vec\gy\nfc [\nab^\fI_n\vec\ga]$. Since $v(\nab_n^\fI\vec\ga)=C_{\nab_n^\fI\vec\ga}$ we have $\vec\gy\notin C_{\nab_n^\fI\vec\ga}$, hence there is an $m\leq n$ such that $\gy_m<\ga_m$. Consider a word $A\in \Wo_m$ such that $o_m(A)=\gy_m$. Recall that the point on the main axis corresponding to $A$ is $\iota(A)=(\gw_m(\gy_m),\dots,\gw^{\gy_m},\gy_m,\ell(\gy_m),\dots)$.

We claim that $\Ig,\vec\gy\nfc \Diamond_m A$, whereas $\Ig,\vec\ga\fc \Diamond_m A$. The former holds, since for all $\vec\gd$ such that $\vec\gy R_m \vec\gd$ one has $\gd_m<\gy_m$, hence $\vec\gd\nleqslant_\fI \iota(A)$ and $\Ig,\vec\gd\nfc A$. On the other hand, $\Ig,\vec\ga\fc \Diamond_m A$ holds, since there is a sequence $\vec\ga':=(\ga_0,\ldots,\ga_{m-1},\gy_m,\gy_{m+1},\dots)$ such that $\vec\ga R_m \vec\ga'$ and $\Ig,\vec\ga'\fc A$.

To show that $\vec\ga'\leq_\fI \iota(A)$ we prove that $\al{i\leq m}\gw_{m-i}(\gy_m)\leq \ga_i$ by downward induction on $i\leq m$. Assume the claim holds for some $i$ such that $0<i\leq m$. Then $\ga_{i-1}\geq \gw^{\ell(\ga_{i-1})}\geq \gw^{\ga_i}\geq \gw^{\gy_i}=\gy_{i-1}$.
\ep

In order to define the operations $\Diamond^{\Rc}_n$ on the set of bounded $\Rc$-theories we need a few definitions. An $\Rc$-theory $T$ is of \emph{level $n$} if $T$ is generated by a (nonempty) set of formulas $\Diamond_n A$ such that $A\in \Wo_n$. A theory $T$ is \emph{of level at least $n$} if it is generated by a (nonempty) subset of $\Wo_n\setminus\{\top\}$.

\bl \label{strat} Every bounded $\Rc$-theory $T$ is representable in the form $T=T_0\land_\Rc T_1\land_\Rc\dots\land_\Rc T_n$ where each $T_i$ is of level $i$ or $T_i=\top_\Rc$.
\el

\bp\ Recall that every $\Rc$-formula is $\Rc$-equivalent to an ordered formula. Moreover, every variable-free $\Rc$-formula in which only the modalities $\Diamond_i$ with $i\geq m$ occur is equivalent to a word in $\Wo_m$. Hence, every formula is equivalent to a conjunction of formulas of the form $\Diamond_i A$ with $A\in\Wo_i$. Since $T$ is bounded, the set of indices of modalities occurring in the axioms of $T$ is bounded, say by $n$. Hence, each axiom of $T$ can be replaced by a finite set of formulas of various levels below $n$ and one can partition the union of all these axioms into the disjoint subsets of the same level. \ep

\bl \label{wn} For each $\vec\ga\in \mI$ such that $\ga_n>0$, the theory generated by $[\vec\ga]\cap \Wo_n$ corresponds to the sequence $\vec\ga':=(\gw_n(\ga_n),\dots,\gw^{\ga_n},\ga_n,\ga_{n+1},\dots)$.
\el

We remark that if $\ga_n=0$ then the theory generated by $[\vec\ga]\cap \Wo_n$ is $\top_\Rc$.
\bp\ Let $T$ be the theory generated by $[\vec\ga]\cap \Wo_n$. We consider a $\vec\gb\in \mI$ such that $[\vec\gb]=T$ and show that $\vec\gb=\vec\ga'$.
It is easy to see that $\vec\ga\leq_\fI\vec\ga'$ and that the submodel of $\Ig$ generated from $\vec\ga$ by the relations $R_k$, for all $k\geq n$, is isomorphic to the submodel generated by these relations from $\vec\ga'$. Hence, if $B$ is a formula in which only the modalities $\Diamond_k$ with $k\geq n$ occur, then $\Ig,\vec\ga\fc B$ holds iff $\Ig,\vec\ga'\fc B$. It follows that $[\vec\ga]\cap \Wo_n\subseteq [\vec\ga']$, that is, $\vec\ga'\leq_\fI\vec\gb$.

Now assume $\vec\ga'<_\fI \vec\gb$, so there is a $k\in\gw$ such that $\gb_k<\ga'_k$. If $k<n$ then $\gb_k<\gw_{n-k}(\ga_n)$. For all ordinals $\gy,\gd$, if $\gy<\gw^\gd$ then $\ell(\gy)<\gd$. Then, by induction, for all $i=k,\dots,n$ we obtain $\gb_{i}<\gw_{n-i}(\ga_n)$. Ergo $\gb_n<\ga_n$.

So, we may assume that $k\geq n$. In this case consider a word $B\in \Wo_k$ such that $o_k(B)=\gb_k+1$. Then,
$$\iota(B)=(\gw_k(\gb_k+1),\dots,\gw_1(\gb_k+1),\gb_k+1,0,\dots).$$
We have $\Ig,\vec\gb\nfc B$, since $\gb_{k}+1>\gb_k$. On the other hand, $$\al{i\leq k} \gw_i(\gb_k+1)\leq \ga_{k-i},$$ which is easy to see by induction on $i$. It follows that $\Ig,\vec\ga\fc B$, therefore $[\vec\gb]\neq T$, a contradiction.\ep

\bcor
For each $\vec\ga\in \mI$, $[\vec\ga]$ is of level at least $n$ iff $\ga_n>0$ and \beq \al{i<n} \ga_i=\gw_{n-i}(\ga_n). \label{om} \eeq
\ecor

\bl\
For each bounded $\Rc$-theory $T$ of level at least $n$, there is an $\Rc$-formula $A\in\Wo_n$ such that $\nab^{\Rc}_n A=\nab^{\Rc}_n T$ holds in $\fT^0_\Rc$.
\el

\bp\ Suppose $T=[\vec\ga]$ is of level at least $n$.
Let $A\in\Wo_n$ be such that $o_n(A)=\ga_n>0$. Then, by Lemma \ref{iso-nab-th},   $\nab^{\Rc}_n(T)=\nab^{\Rc}_n([\vec\ga])=[\nab^\fI_n\vec\ga]$. By \refeq{om} we have  $$\nab^\fI_n\vec\ga=(\gw_n(\ga_n),\gw_{n-1}(\ga_n),\dots,\ga_n,0,\dots).$$
On the other hand, $\iota(A)=(\gw_n(\ga_n),\gw_{n-1}(\ga_n),\dots,\ga_n,\ell(\ga_n),\dots),$ and we obtain
$\nab^{\Rc}_n A = [\nab^\fI_n(\iota(A))]=[(\gw_n(\ga_n),\gw_{n-1}(\ga_n),\dots,\ga_n,0,\dots)]$.
Proposition \ref{iso-th} yields the result. \ep

Now we can give the following definition of the theory $\Diamond^{\Rc}_n T$, for each bounded $\Rc$-theory $T$.

If $T$ is of level at least $n$ or $T=\top_\Rc$, we let
$\Diamond^{\Rc}_n T$ be the theory generated by the formula $\Diamond_n A$, where $A\in\Wo_n$ is such that $\nab^{\Rc}_n A=\nab^{\Rc}_n T$ in $\fT_\Rc^0$. (Notice that this definition is correct, since any two words $A_1,A_2$ satisfying $\nab^{\Rc}_n A_1=\nab^{\Rc}_n A_2$ in $\fT_\Rc^0$ also satisfy $\Diamond_n A_1=_\Rc\Diamond_n A_2$ by Corollary \ref{nab-inj}.)

For each $i\leq n$, let $T_i$ denote the theory generated by $T\cap \Wo_i$.
We define
$$\Diamond^{\Rc}_n(T):=\Diamond^{\Rc}_0(T_0)\land_\Rc \Diamond^{\Rc}_1(T_1)\land_\Rc \dots \land_\Rc \Diamond^{\Rc}_n(T_n).$$

The following lemma shows that this definition agrees with the operations on the Ignatiev algebra.

\bl For all $\vec\ga\in\fI$, $\Diamond^{\Rc}_n([\vec\ga])=[\Diamond^\fI_n(\vec\ga)]$.
\el

\bp\ If $T=[\vec\ga]$ then by Lemma \ref{wn}, for each $i\leq n$, either the theory $T_i:=T\cap \Wo_i$ is $\top_\Rc$ or corresponds to the sequence $\vec\ga':=(\gw_i(\ga_i),\dots,\gw^{\ga_i},\ga_i,\ga_{i+1},\dots)$ with $\ga_i>0$. If $T_i=\top_\Rc$ we have $\Diamond_i^\Rc T_i=\Diamond_i\top$. Otherwise, $\Diamond^{\Rc}_iT_i=\Diamond_i A_i$ where $A_i$ corresponds to $(\gw_i(\ga_i),\dots,\gw^{\ga_i},\ga_i,\ell(\ga_i),\dots)$. In both cases  $$\Diamond^{\Rc}_i T_i= [(\gw_i(\ga_i+1),\dots, \gw^{\ga_i+1},\ga_i+1,0,\dots)].$$
Then we observe that $\Diamond^{\Rc}_n(T)=\Diamond^{\Rc}_0(T_0)\land_\Rc \Diamond^{\Rc}_1(T_1)\land_\Rc \dots \land_\Rc \Diamond^{\Rc}_n(T_n)$ corresponds to the cone generated by $(\ga_0+1, \ga_1+1,\dots,\ga_n+1,0,\dots)$ in $\cE$ which coincides with the cone of $\Diamond_n^\fI(\vec\ga)$ (cf.\ Lemma \ref{Diam-hom}).
\ep

Using Lemma~\ref{Diam-hom} we can also isomorphically represent $\fI$ as an algebra of cones in $\Ig$. Given a cone $C\in \mC(\Ig)$ let $\Diamond^\fC_n(C):=\bigcap_{i\leq n} R_i^{-1}(C)$. We also define $$\nab_n^\fC(C):=\bigcap\{R_i^{-1}(D):D\in\fC(\Ig),\ i\leq n,\ R_i^{-1}(D)\supseteq C\}.$$

We summarize the main results of this paper in the following theorem.
\bt \label{iso}
The following structures are isomorphic:
\benr
\item $\ofG^0_T$, for any sound G\"odelian extension $T$ of $\EA$;
\item $\fL^0_{\Rcn}$, the Lindenbaum--Tarski algebra of the variable-free fragment of $\Rcn$;
\item $\fI=(\mI,\land_\fI,\{\Diamond^\fI_n,\nab^\fI_n:n\in\gw\})$;
\item $(\fT^0_\Rc,\land_\Rc,\{\Diamond^{\Rc}_n,\nab^{\Rc}_n:n\in\gw\})$;
\item $\fC(\Ig)=(\mC(\Ig),\cap, \{\Diamond^\fC_n,\nab_n^\fC:n\in\gw\})$.
\eenr
\et

\bp\ We only need to prove the isomorphism of (v) with either (iii) or (iv). Proposition~\ref{iso-th} provides the isomorphisms of the semilattice reducts.
Further, for all $\vec\ga\in \mI$, $\Diamond^\fC_n(C_{\vec\ga})=C_{\Diamond^\fI_n(\vec\ga)}$
by Lemma~\ref{Diam-hom} (i). Hence, $\Diamond^\fC_n$ corresponds to $\Diamond_n^\fI$ of (iii). On the other hand, $\nab_n^\fC(C_{\vec\ga})=v(\nab_n^\Rc([\vec\ga]))$. Hence, $\nab_n^\fC$ corresponds to $\nab_n^\Rc$ of (iv).
\ep

We remark that the algebra $\fC(\Ig)$ has rather simple definitions of meet and diamonds, but somewhat convoluted nablas. In contrast, $\fT^0_\Rc$ has simple meet and nablas but somewhat convoluted diamonds. The algebra $\fI$, perhaps the most elegant of all three, has a more complicated meet operation (though the order relation $\leq_\fI$ is simple).

Finally, we briefly return to the subject of conservativity spectra and look at it from the point of view of established isomorphisms.

Let us call a theory $S$ in the language of $\PA$ \emph{bounded} if $S$ is contained in a consistent finitely axiomatizable theory. The unboundedness theorem by Kreisel and L\'evy \cite{KrL} yields that $\ord_n(S)=0$, for all sufficiently large $n\in\gw$, whenever $S$ is bounded. We need to restrict ourselves to bounded subtheories of $\PA$ if we want to establish a bijection between their conservativity spectra and the Ignatiev algebra.

\bt\ \benr \item Let $T$ be a G\"odelian extension of $\EA^+$ and let $\vec\ga$ be the conservativity spectrum of $T$. If $\PA\vdash T$ then $\vec\ga\in \oI$. If, in addition, $T$ is bounded, then $\vec\ga\in I$.
\item Let $\vec\ga\in \fI$, $A$ be a variable-free $\Rcn$-formula corresponding to $\vec\ga$ via the isomorphism, and $A^*\in\fG^0_{\EA^+}$ its arithmetical interpretation. Then $A^*$ is a bounded subtheory of $\PA$ and $\vec\ga$ is the conservativity spectrum of $A^*$.
\item Under the same assumptions, $A^*$ is the weakest theory with the given conservativity spectrum $\vec \ga$.
\eenr
\et

\bp\ (i) In view of Lemma \ref{spec-ign}, for the first claim it is sufficient to prove that $\al{n\in\gw} \ga_n\leq\ge_0$. Since $\PA$ contains $T$, this follows from Proposition~\ref{insense}~(i).

Since $\PA$ is equivalent to the union of theories $\{\mR_n(1): n\in\gw\}$, any finite subtheory of $\PA$ is contained in a theory of the form $\mR_n(1)$, for some $n\in\gw$ (we write $1$ for $1_{\EA^+}$). Hence, its conservativity spectrum is $\leq_\fI$ above that of $\mR_n(1)$, that is, belongs to $\fI$.

(ii) That $A^*$ is a bounded subtheory of $\PA$ easily follows by induction on the build-up of $A$. The equality $\spec(A^*)=\vec \ga$ is a part of Theorem~\ref{spect}.

(iii) This follows from the fact that any theory $T$ such that $\spec(T)\leq_\fI \vec\ga$ must contain the union of progressions \[\mR_0^{\ga_0}(1)\land \mR_1^{\ga_2}(1)\land\dots\land \mR_k^{\ga_k}(1),\]
which is equivalent to $A^*$ by Theorem~\ref{spect}.
\ep

Let $\spec(T)$ denote the conservativity spectrum of $T$ and let $\thh:\fI\to \ofG_{\EA^+}$ denote the natural isomorphic embedding of $\Rcn$-algebras. As we already noted, $\thh(\vec\ga)$ is a bounded subtheory of $\PA$, for each $\vec\ga$.

\bcor The maps $\thh$ and $\spec$ form a Galois connection: for each bounded subtheory $S$ of $\PA$,
$$\spec(S)\leq_\fI \vec\ga \iff S\leq_{\EA^+} \thh(\vec\ga).$$
\ecor

We remark that the map $\spec$ is order-preserving, however it is not a semilattice homomorphism, even when restricted to bounded subtheories of $\PA$. For example, it is well known that $\ord_1(I\Sigma_1)=\gw=\ord_1(I\Pi_2^-)$ and both theories are $\Pi_2^0$-regular:
\begin{eqnarray*}\spec(I\Pi_2^-)& = & (\gw^\gw,\gw,0,\dots) \\
\spec(I\Sigma_1) & = & (\gw^\gw,\gw,1,0,\dots).\end{eqnarray*}
On the other hand,  $\ord_1(I\Sigma_1\land_\EA I\Pi_2^-)=\gw^2>\gw$ and
$$\spec(I\Sigma_1\land_\EA I\Pi_2^-)=(\gw^{\gw^2},\gw^2,1,0,\dots).$$

\section{A universal Kripke frame for the variable-free fragment of $\Rcn$} \label{ukf}

In view of Theorem \ref{iso} it is natural to ask if one can describe a convenient universal Kripke frame for the variable-free fragment of $\Rcn$. There are two known general constructions associating with an SLO $\fB=(B,\land^\fB,\{a^\fB:a\in\Sigma\})$ its `dual' Kripke frame, so that $\fB$ is embeddable into the algebra of subsets of that frame (see \cite[Section 4.1]{KKTWZ}). One construction is similar to the way the canonical model of a strictly positive logic $L$ is obtained from its Lindenbaum--Tarski algebra and goes from $\fB$ to the set of all filters of $\fB$ equipped with binary relations $\{R_a:a\in \Sigma\}$ such that, for all filters $F,G$,
$$
F R_a G \iffdef \al{x\in G} a^\fB(x)\in F.
$$
The corresponding frame for the $\Rcn$-algebra $\fI$ is constructively described in \cite{Bek18c} in terms of appropriate sequences of ordinals. However, the relations of the frame look sufficiently complicated, so that one would really want a simpler construction for practical use.

Another approach (see~\cite{Jack,KKTWZ}) is to consider the set $B$ itself as a dual space, and to specify binary relations on $B$ by
$$
x R_a y \iffdef x\leq_\fB a^\fB (y).
$$
Let $\fB^*$ denote the Kripke frame $(B,\{R_a:a\in\Sigma\})$ together with the canonical valuation $v:\fB\to \cP(B)$, where $v(x):=\{y\in B:y\leq_\fB x\}$.

\bl
For all $x,y\in \fB$ and $a\in\gS$, the following relations hold in $\fB^*$:
\benr
\item $v(x\land_\fB y)=v(x)\cap v(y)$; \item $R_a^{-1}(v(x))=v(a^{\fB}(x))$.
\eenr
\el

\bp\ Claim (i) is just the fact that $x\land_\fB y$ is the g.l.b.\ of $x$ and $y$. To prove Claim (ii) we argue as follows: $z\in  R_a^{-1}(v(x))$ means there is a $u\leq_\fB x$ such that $z R_a u$, that is, $z\leq_\fB a^{\fB}(u)$. Thus, if $z\in  R_a^{-1}(v(x))$ we have by monotonicity $a^{\fB}(u)\leq_\fB a^{\fB}(x)$ and therefore $z\leq_\fB a^{\fB}(x)$.

If $z\leq_\fB a^{\fB}(x)$ then we take $x$ for $u$ and observe that $u\leq_\fB x$ and $z\leq_\fB a^{\fB}(u)$, hence $z\in  R_a^{-1}(v(x))$. \ep

We obtain the following corollaries.

\bpr \benr \item The map $v:\fB\to \cP(B)$ is an embedding of $\fB$ into the algebra $(\cP(B),\cap,\{R^{-1}_a:a\in\Sigma\})$.
\item If $A,B$ in $\cL_\gS$ are variable-free, then $A\vdash B$ holds in $\fB$ iff $\fB^*,x\fc A\to B$ for all $x\in B$.
    \eenr
\epr

\bcor The variable-free fragment of $\Rcn$ is complete w.r.t.\ $\fI^*$. \ecor

The Kripke frame $\fI^*$ has a simple constructive characterization. We know that its domain is the set $I$ of all sequences of ordinals $\vec\ga=(\ga_0,\ga_1,\dots)$ such that, for all $n\in\gw$, $\ga_n<\ge_0$ and $\ga_{n+1}\leq \ell(\ga_n)$. Our task is to characterize the relations $R^*_n$ and $S^*_n$ on $I$ corresponding to, respectively, $\Diamond_n$ and $\nab_n$, for all $n\in\gw$, where
\begin{eqnarray*}
\vec\ga R^*_n \vec\gb & \iff & \vec\ga\leq_\fI \Diamond^\fI_n\vec\gb; \\
\vec\ga S^*_n \vec\gb & \iff & \vec\ga\leq_\fI \nab_n^\fI\vec\gb.
\end{eqnarray*}
The answer is given by the following proposition.

\bpr For all $\vec\ga,\vec\gb\in I$,
\benr
\item $\vec\ga R^*_n \vec\gb \iff \al{i\leq n} \ga_i>\gb_i;$
\item $\vec\ga S^*_n \vec\gb \iff \al{i\leq n} \ga_i\geq\gb_i.$
\eenr
\epr

\bp\ Claim (ii) is obvious, since $\nab^\fI_n\vec\gb=(\gb_0,\gb_1,\dots,\gb_n,0,\dots)$. To prove Claim (i) we recall that $\Diamond^\fI_n\vec\gb=(\gb'_0,\gb'_1,\dots,\gb'_n,0,\dots)$ where
$\gb'_{i}=0$ for $i>n$ and $\gb'_i=\gb_i+\gw^{\gb'_{i+1}}$ for $i\leq n$. Clearly, for all $i\leq n$ $\gb'_i>\gb_i$. Hence, the `only if' part of the claim is obvious.

To prove the `if' part, we assume $\al{i\geq n}\ga_i> \gb_i$ and prove by downwards induction on $i\leq n$ that $\al{i\geq n}\ga_i\geq \gb'_i$. If $i=n$ then $\gb'_i=\gb_i+1$ and the claim is clear. If $i<n$ then $\ga_i>\gb_i$ and by the induction hypothesis $\ga_{i+1}\geq \gb'_{i+1}$. Since $\vec\ga\in I$ we have $\ell(\ga_i)\geq \ga_{i+1}\geq \gb'_{i+1}$. At this point we need an auxiliary lemma.

\bl\ For any ordinals $\ga,\gb,\gy$, if $\ga>\gb$ and $\ell(\ga)\geq \gy$ then $\ga\geq \gb+\gw^\gy$.
\el
\bp\ We can write $\ga=\gb+\nu$ with $\nu>0$. Then $\ell(\nu)=\ell(\ga)\geq\gy$, hence $\nu\geq \gw^\gy$ and $\ga=\gb+\nu\geq \gb+\gw^\gy$. \ep

By this lemma we conclude that $\ga_i\geq \gb_i+\gw^{\gb'_{i+1}}=\gb'_i$ and the induction step is complete. \ep

Looking at the frame $\fI^*$ as a dual of the Lindenbaum--Tarski algebra of the variable-free fragment of $\Rcn$ we observe that, for any $A,B\in \Fo^\nab$, $A R^*_n B$ holds iff $A\vdash_\Rcn \Diamond_n B$. Hence, $R_n^*$ is the same as the previously considered relation $<_n$ on words (now extended to all variable-free formulas of $\Rcn$).

On the other hand, $A S^*_n B$ holds iff $A\vdash_\Rcn \nab_n B$. Hence, $S^*_n$ is the same as the $\Pi^0_{n+1}$-conservativity relation previously denoted $\vdash_n$ (cf Section \ref{var-free}).

Instead of the Lindenbaum--Tarski algebra of the variable-free fragment of $\Rcn$ we can also work directly with its isomorphic arithmetical counterpart, the SLO $\ofG_T^0$ for a sound extension $T$ of $\EA$. Then, $\gs R_n^* \nu$ means that the G\"odelian theory $\gs$ proves $\mR_n(\nu)$, and $\gs S_n^* \nu$ means that $\nu$ is $\Pi_{n+1}^0$-conservative over $\gs$.

\brem The same definitions also apply to a much larger Kripke frame $\ofG^*_T$ that is dual to the $\Rcn$-algebra of all G\"odelian extensions of $T$, $\ofG_T$.
\erem

\brem
A recent paper by Hermo Reyes and Joosten~\cite{JooRey17} introduces a universal Kripke frame for the so-called Turing--Schmerl Calculus. This model turns out to be very similar to $\fI^*$. The differences amount to the following two aspects. Firstly, their relations $R_n$ can be defined as $R^*_n\cap \leq_\fI$. This reflects the fact that all their modalities satisfy the principle $\Diamond A\vdash A$. Secondly, their models lack the $S_i$ relations, but allow the $\alpha$-iterations of relations $R_n$.
\erem

\begin{appendices}

\section{Irreflexivity of $<_0$ in $\Rc$} \label{irreflex}

We work in (the variable-free fragment of) the reflection calculus $\Rc$. We will use the techniques of Kripke models for $\Rc$. The notions of the \emph{canonical tree} for a formula $A$, its \emph{RC-closure} $\Rc[A]$ and that of an \emph{RC-model} are defined in \cite{Bek14}. We recall that $\Rc[A]$ is an RC-model satisfying $A$ at the root. Its valuation will be empty if $A$ is variable-free.

The following lemma is easily obtained from Lemma~\ref{ord} taking into account that words in $\Wo_n$ are linearly pre-ordered by $<_n$.

\bl Any variable-free formula of $\Rc$ is equivalent to $\top$ or to a
formula of the form $A\circeq\bigwedge_{i\leq k} \Diamond_{m_i} A_i$ where

\benr \item $A_i\in \Fo_{m_i}$, for each $i$;
\item $m_0 > m_1 >  \dots > m_k$;
\item $\Diamond_{m_i} A_i\nvdash_\Rc \Diamond_{m_j} A_j$, for all $j>i$.
\eenr \el

Such formulas are called \emph{properly ordered}. If $A$ is
properly ordered, then $\Rc[A]$ can be characterized as follows.

If $A\circeq\top$ then $\Rc[A]$ is the irreflexive singleton frame. If
$A\circeq\bigwedge_{i\leq k} \Diamond_{m_i} A_i$ then $\Rc[A]$ consists of the
disjoint union of the frames $\Rc[A_i]$, for all $i\leq k$,
augmented by a new root $a$. In addition to all the relations
inherited from the frames $\Rc[A_i]$, the following relations are
postulated:

\ben
\item $a R_{n} x$, for each $i\leq k$, $n\leq m_i$ and $x\in \Rc[A_i]$;
\item $x R_n y$, for each $i\leq k$, $n<m_i$, and $x,y\in
\bigcup_{j\leq i} \Rc[A_j]$;
\item $x R_n y$, for each $i\leq k$, $n\leq m_i$, $y\in \Rc[A_i]$
and $x\in \bigcup_{j<i} \Rc[A_j]$. \een

The following lemma is routine.
\bl $\Rc[A]$ thus described is an RC-frame. \el

\bt For any formula $A$ of $\Rc$, $A\nvdash_{\Rc}\Diamond_0 A$. \et

\bp\ It is sufficient to prove the claim for variable-free and properly ordered $A$. For such an $A$, we argue by induction on the length of $A$. Basis is trivial. Suppose $A=\bigwedge_{i\leq k} \Diamond_{m_i} A_i$. If $A\vdash \Diamond_0A$
then there is a homomorphism $f$ of $\Rc[A]$ into itself such that
$a R_0 f(a)$. Then there is an $i\leq k$ such that $f(a)\in
\Rc[A_i]$.

Let $X$ denote the subset of $\Rc[A]$ corresponding to $\Rc[A_i]$. Consider any $n\geq
m_i$ and an $R_n$-arrow whose source is in $X$. By the construction
of $\Rc[A]$, this arrow can only be an old arrow from the frame
$\Rc[A_i]$. Hence, the target of this arrow will also be in $X$.
Since $A_i\in \Fo_{m_i}$, it follows that $f(X\cup\{a\})\subseteq
X$. The subset $X\cup\{a\}$ together with all the inherited
relations can be considered as a submodel of $\Rc[A]$ isomorphic to
$\Rc[\Diamond_{m_i}A_i]$. Hence, $f$ induces a homomorphism $f: \Rc[\Diamond_{m_i}A_i]\to
\Rc[A_i]$. This implies that either $A_i\vdash_{\Rc} \Diamond_{m_i} A_i$ (if
$f(a)$ is the root of $\Rc[A_i]$), or $A_i\vdash_{\Rc}
\Diamond_{m_i}\Diamond_{m_i}A_i\vdash_{\Rc} \Diamond_{m_i}A_i$ (if $f(a)$ is strictly above the root). In any case $A_i\vdash_{\Rc} \Diamond_{m_i}A_i\vdash_{\Rc} \Diamond_0A_i$
contradicting the induction hypothesis. \ep

\section{Uniform definability of computable operators}

\bt  An operator $R:\fG_\EA\to \fG_\EA$ is uniformly definable iff $R$ is computable. \et

\bp\ The main point is to show that computable $R$ are uniformly definable. Let $R$ be computable, hence there is a $\Sigma_1^0$-formula $\Ax_R(x,y)$ such that $\Ax_R(x,\ol{\gn{\gs}})$ numerates the theory  $R(\gs)$ for each $\gs$. Notice that $R(\gs)$ is an elementary formula, for each $\gs$. We claim that one can select $\Ax_R$ in such a way that for each $\gs$ there is an elementary numeration $\gd$ such that \beq \EA\vdash \al{x}(\Ax_R(x,\ol{\gn{\gs}})\eqv \gd(x)).\label{gd} \eeq
Let $\Sat_{\Delta_0}(e,x)$ be a $\Sigma_1^0$-truthdefinition for elementary formulas that can be represented in the form
$$\Sat_{\Delta_0}(e,x) \eqv \ex{q\leq 2_{d(e)}^x} \ T(e,x,q),$$
where $T(e,x,q)$ is an elementary formula expressing that $q$ is a protocol of a computation verifying that an elementary formula $e$ holds on assignment $x$. For each specific formula $e$, the size of $q$ is bounded by a $d$-fold iterate of exponential function in $x$ where $d$ elementarily depends on $e$. Whereas in $\EA$ one cannot prove that $2_{d(e)}^x$ is defined for all $e$ and $x$, it is known that for each specific $n$ there is an $\EA$-proof of $\al{x} \ex{y} \: 2_{\bar n}^x=y$. So, for each specific formula $\gs$ there is a number $n=d(\gn{\gs})$ such that provably in $\EA$  \beq \al{x}(\Sat_{\Delta_0}(\ol{\gn{\gs}},x)\eqv  \ex{q}\leq 2_{\bar n}^x\ T(\ol{\gn{\gs}},x,q)). \label{sat} \eeq
Now, if $F_R(x,y)$ is a $\Sigma_1^0$-formula strongly representing the map $R:\gn{\gs}\mapsto \gn{R(\gs)}$, we can define
$$\Ax_R(x,y)\iffdef \ex{e}(F_R(y,e)\land \Sat_{\Delta_0}(e,x)).$$ Then, for each $\gs$ there is a provably unique $\tau=R(\gs)$ such that $\EA\vdash F_R(\ol{\gn{\gs}},\ol{\gn{\tau}})$. Hence,
$\Ax_R(x,\ol{\gn{\gs}})$ is provably equivalent to $\Sat_{\Delta_0}(\ol{\gn{\tau}},x)$ which is equivalent to an elementary formula by \refeq{sat}. This proves \refeq{gd}.

To provide a uniform definition of $R$ we apply a version of Craig's trick and let
$$\Ax_{R'}(x,y)\iffdef \ex{z,p\leq x}(x=\mathrm{disj}(z,\gn{\ol{p}\neq\ol{p}})\land W_R(z,y,p)),$$ where $W_R(z,y,p)$ is an elementary formula expressing that $p$ witnesses $\Ax_R(z,y)$. Here, we may assume that $\EA\vdash W_R(z,y,p)\to z\leq p$.
Clearly, $\Ax_{R'}(x,y)$ is elementary and condition (ii) is satisfied. Externally, it numerates the same family of theories as $\Ax_{R}(x,y)$.
We show that, for each $\gs$,
$$\EA\vdash \al{x} (\Box_{R(\gs)}(x) \to \Box_{R'(\gs)}(x)).$$
First, we obtain an elementary numeration $\gd$ such that
$\EA\vdash \al{x}(\Ax_R(x,\ol{\gn{\gs}})\eqv \gd(x)).$
It follows that $\EA\vdash \al{x} (\Box_{R(\gs)}(x) \eqv \Box_{\gd}(x)).$
Thus, using $\Pi_2^0$-conservativity of $\BS_1$ over $\EA$ it is sufficient to prove $$\EA+\BS_1\vdash \al{x} (\Box_{\gd}(x)\to \Box_{R'(\gs)}(x)).$$
Using $\BS_1$ it is sufficient to prove that $\EA\vdash \al{x}(\gd(x)\to \Box_{R'(\gs)}(x)).$ Reason in $\EA$: Assume $\gd(x)$ then $\Ax_R(x,\ol{\gn{\gs}})$. Hence, there is a witness $p$ such that $W_R(x,\ol{\gn{\gs}},p)$. Then for $u:=\text{disj}(x,\gn{\ol{p}\neq\ol{p}})$ we have $\Ax_{R'}(u,\ol{\gn{\gs}})$ and from $p$ we obtain a proof of $\ol{p}\neq\ol{p}$ and hence a proof of $x$ from hypothesis $u$ in an elementary way. Therefore, $\Box_{R'(\gs)}(x)$. \ep

\end{appendices}

\bibliographystyle{plain}
\bibliography{ref-all2}

\end{document}